\documentclass[11pt,a4paper]{article}
\usepackage[margin=1in]{geometry}
\usepackage[T1]{fontenc}
\usepackage{lmodern}
\usepackage{amsmath,amssymb,amsthm,mathtools}
\usepackage{booktabs,array,enumitem}
\usepackage{microtype}
\usepackage{xcolor}
\usepackage[colorlinks=true,linkcolor=blue!45!black,citecolor=blue!45!black,urlcolor=blue!45!black]{hyperref}
\hypersetup{
 pdftitle={Threshold Representations and Realization Fibers of Specialization in Chang Fuzzy Topologies},
 pdfauthor={S. M. Elsayed},
 pdfsubject={Specialization relations, threshold matrices, and realization fibers in Chang fuzzy topologies}
}
\newtheorem{theorem}{Theorem}[section]
\newtheorem{proposition}[theorem]{Proposition}
\newtheorem{corollary}[theorem]{Corollary}
\newtheorem{lemma}[theorem]{Lemma}
\theoremstyle{definition}
\newtheorem{definition}[theorem]{Definition}
\newtheorem{example}[theorem]{Example}
\theoremstyle{remark}
\newtheorem{remark}[theorem]{Remark}
\newcommand{\F}{\mathcal F}
\newcommand{\FP}{\operatorname{FP}}
\newcommand{\Pp}{\mathsf P}
\newcommand{\Qq}{\mathsf Q}
\newcommand{\K}{\mathsf K}
\newcommand{\cl}{\operatorname{cl}}
\newcommand{\Fix}{\operatorname{Fix}}
\newcommand{\preM}{\preccurlyeq_{m}}
\newcommand{\preQ}{\preccurlyeq_{q}}
\newcommand{\tauM}{\tau_m}
\newcommand{\tauQ}{\tau_q}
\newcommand{\tauC}{\tau_{\Gamma}}
\newcommand{\Up}{\operatorname{Up}}
\title{\textbf{Threshold Representations and Realization Fibers of Specialization in Chang Fuzzy Topologies}}
\author{S. M. Elsayed\\
\small Department of Mathematics, Faculty of Science, Aswan University, Aswan, Egypt\\
\small \texttt{profsalah55@aswu.edu.eg}}
\date{}
\begin{document}
\maketitle

\begin{abstract}
We investigate the information retained by fuzzy-point specialization in Chang fuzzy topologies and the fibers realizing a prescribed G\"odel specialization relation. Every reflexive, minimum-transitive fuzzy preorder has a greatest stratified Chang realization. Its realization fiber is a singleton precisely for the constant relation $1$; otherwise, below each member lies a continuum-sized descending chain with no lower bound in the fiber, also within the Alexandrov subclass. On finite support sets, every nontrivial stratified fiber has cardinality continuum without assuming Alexandrov completeness; this remains true for Alexandrov fibers on countable sets. We describe the entire two-point Alexandrov fiber order as the reverse of a product of classes of join-preserving threshold functions. For arbitrary Chang topologies, quasi-specialization determines membership specialization and its kernel completion, which in turn determine G\"odel specialization; the converse determinations fail in general. A central reconstruction result identifies the greatest quasi-specialization realization by \emph{exact agreement on every one- or two-point restriction}, equivalently by closure in a product of discrete grade intervals. In particular, quasi-specialization determines every Chang topology on a finite support set. We also represent general Alexandrov Chang topologies by matrices of join-preserving threshold maps and isolate the additional scalar compatibility law required for reconstruction from the graded relation. These results separate information encoded by finite-coordinate point comparisons from global restrictions on fuzzy open families.
\end{abstract}

\section{Introduction}
A Chang fuzzy topology is a family of fuzzy subsets closed under arbitrary joins and finite meets \cite{Chang1968}. Membership and quasi-coincidence represent that family on the ordinary point set $E=X\times(0,1]$. Point spaces and attachment representations are well established \cite{PuLiu1980,Eroglu1990,Eroglu1991,Guido2010,GuidoScarciglia2011}. With fixed point labels and grade coordinates, either point topology determines the original fuzzy open family. Passing to specialization preorders, however, may discard information.

We study three levels of order information: the quasi-point preorder, the membership-point preorder (which compares possibly unequal grades), and the G\"odel specialization relation $r_\Omega(x,y)$, which assigns a single grade to each ordered pair of supports. Our tool is the open-kernel operator $K_\Omega$, related to the kernel completions studied in \cite{ElsayedJCAM2027}, and its threshold profiles
\[
 h_{xy}(t)=K_\Omega(x_t)(y),\qquad x,y\in X,\quad t\in[0,1].
\]
The profiles encode precisely the membership preorder. The inverse problem is to recover them, or the original topology, from the less detailed graded relation, keeping the underlying labeled set fixed.

The main contributions are threefold.
\begin{enumerate}[label=\textup{(\roman*)},leftmargin=*]
\item \emph{General realization fibers.} Theorem~\ref{thm:generalfibre} establishes a dichotomy for all stratified Chang realizations of a prescribed G\"odel fuzzy preorder: a unique constant topology, or continuum-sized descending chains with no lower bound strictly below every realization. The construction preserves the Alexandrov property, and the cardinality of nontrivial Alexandrov fibers is exact on at most countable sets.
\item \emph{Threshold representations.} Theorem~\ref{thm:matrixrep} represents general Alexandrov Chang topologies by matrices of join-preserving maps. Theorem~\ref{thm:firstdeparture} identifies what G\"odel specialization retains from a stratified threshold profile, and Theorem~\ref{thm:fibreorder} describes the entire two-point Alexandrov fiber order.
\item \emph{Reconstruction and cardinality.} Theorem~\ref{thm:quasireconstruction} reconstructs the membership kernel from quasi-specialization. Theorem~\ref{thm:finitequasi} establishes exact reconstruction on finite support sets, and Theorem~\ref{thm:qmax} identifies the greatest quasi-specialization realization by finite restrictions, equivalently by restrictions to at most two points. Theorem~\ref{thm:finiteallcard} counts all stratified realizations on finite sets. The examples in Section~\ref{sec:families} show that the reconstruction hierarchy is strict in general.
\end{enumerate}

\noindent\textbf{Relation to previous work and novelty boundary.}
Point-space and attachment representations are established
\cite{PuLiu1980,Guido2010,Eroglu1990,Eroglu1991,GuidoScarciglia2011}.
Graded specialization and fuzzy upper-set topologies are established
\cite{LaiZhang2006,Fang2007,YaoShi2008}, as are join-preserving and
approximation operators \cite{Kim2014join,OhKim2021,JarvinenKortelainen2007}.
Fang treats graded openness, whereas our families are ordinary Chang
fuzzy topologies; the intuitionistic-fuzzy kernel theory
\cite{ElsayedJCAM2027} is related. Neither point representation nor
scalar relational kernel formulas are claimed as new. Theorem~\ref{thm:matrixrep}
presents established operator-matrix theory in coordinates suited to
Alexandrov Chang realizations without requiring scalar compatibility.

Topkis \cite[Theorem~1]{Topkis1976} established a two-coordinate
projection principle for sublattices of finite lattice products
(compare the Baker--Pixley majority-term principle). We do not claim
that lattice projection principle as new. Theorem~\ref{thm:qmax}
instead characterizes the greatest Chang realization of a fixed
quasi-point preorder through exact finite-coordinate closure,
equivalently one- and two-coordinate matching, also on infinite
support sets. This relates a classical lattice principle to
quasi-specialization reconstruction.

Lai--Zhang \cite[Theorem~4.8 and Example~4.24]{LaiZhang2006}
give the finest fuzzy-preorder realization and a specific failure
of a coarsest realization. Theorem~\ref{thm:generalfibre} proves a
uniform strengthening: for every nontrivial reflexive,
minimum-transitive relation, below every stratified Chang
realization, including a non-Alexandrov one, lies a strict
continuum-sized descending chain with no lower bound in the fiber.
Theorem~\ref{thm:fibreorder} classifies the entire inclusion order
of each two-point Alexandrov fiber. These are the principal
proposed contributions; no exhaustive priority claim is made.

Benoumhani and Jaballah studied finite fuzzy-topology enumeration
\cite{BenoumhaniJaballah2017,BenoumhaniJaballah2019}. Here the
membership lattice is the entire real unit interval and the
objects counted are fibers with fixed G\"odel specialization.
Lemma~\ref{lem:countquasi} records boundary-attainment information
not captured by numerical kernel profiles. Exact matching of
finite-coordinate restrictions always uses the given support and
grade labels. Sections~\ref{sec:background}--\ref{sec:profiles}
give the framework; Sections~\ref{sec:loss}--\ref{sec:twopoint}
treat scalar compatibility and fibers; Section~\ref{sec:families}
treats reconstruction and cardinality; Section~\ref{sec:order-maps}
treats maps; and Appendix~\ref{sec:other} contains auxiliary
point-space results.

\section{Definitions and known kernel facts}\label{sec:background}
Let $X$ be a nonempty set, and put $\F(X)=[0,1]^X$ with the pointwise order. A \emph{Chang fuzzy topology} is a family $\Omega\subseteq\F(X)$ containing the zero and unit functions, closed under arbitrary pointwise suprema and finite pointwise infima. A topology is \emph{stratified} if it contains every constant fuzzy set, and is \emph{Alexandrov} if it is closed under arbitrary pointwise infima. \emph{Throughout, ``Alexandrov'' means arbitrary-meet closure in the Chang sense only.} We do not build stability under scalar G\"odel implication into this definition, even for stratified topologies. That stronger requirement appears separately in Proposition~\ref{prop:residualboundary}; conflating the two notions would eliminate some of the realization fibers studied here. Following the fuzzy-point convention \cite{PuLiu1980}, the fuzzy point $x_\alpha$ is the function taking value $\alpha$ at $x$ and zero elsewhere, where $0<\alpha\le1$. Write
\[
 E=\FP(X)=\{x_\alpha:x\in X,\ 0<\alpha\le1\}.
\]
For $A\in\F(X)$, let
\[
 \Pp(A)=\{x_\alpha\in E:\alpha\le A(x)\},\qquad
 \Qq(A)=\{x_\alpha\in E: A(x)>1-\alpha\}.
\]
The first construction uses membership, the second quasi-coincidence. For $T\subseteq E$, put
\[
 J(T)(x)=\sup\bigl(\{\alpha:x_\alpha\in T\}\cup\{0\}\bigr).
\]
Then $J(\Pp(A))=A$. The map $\Pp$ preserves arbitrary meets and binary joins, but it need not preserve infinite joins. The map $\Qq$ preserves arbitrary joins and finite meets:
\begin{equation}\label{eq:Q-ops}
 \Qq\Bigl(\bigvee_i A_i\Bigr)=\bigcup_i\Qq(A_i),\qquad
 \Qq(A\wedge B)=\Qq(A)\cap\Qq(B).
\end{equation}
Unlike strict cuts, nonstrict cuts need not preserve arbitrary joins.

We recall the open-kernel operator
\begin{equation}\label{eq:kernel}
 \K(A)=\bigwedge\{U\in\Omega:A\le U\}.
\end{equation}
It is extensive, monotone, idempotent, and preserves arbitrary joins:
\begin{equation}\label{eq:join}
 \K\Bigl(\bigvee_i A_i\Bigr)=\bigvee_i\K(A_i).
\end{equation}
The empty-family case holds because $\K(0)=0$. For completeness, here is a direct justification of the nonempty-family case (the kernel law is also used in \cite{ElsayedJCAM2027}). Given $(A_i)_{i\in I}\subseteq[0,1]^X$, put $\Omega_i=\{U\in\Omega:A_i\le U\}$; each $\Omega_i$ contains $1$. Complete distributivity of $[0,1]^X$ yields
\[
 \bigvee_{i\in I}\bigwedge_{U\in\Omega_i}U
 =\bigwedge_{\varphi\in\prod_{i\in I}\Omega_i}
      \bigvee_{i\in I}\varphi(i).
\]
For each choice $\varphi$, the right-hand inner join belongs to $\Omega$ and dominates $\bigvee_i A_i$, so
$\K(\bigvee_i A_i)\le\bigvee_i\K(A_i)$. The reverse inequality follows from monotonicity of $\K$. Thus \eqref{eq:join} holds. We also use $\K(U)=U$ for $U\in\Omega$. The fixed points of $\K$ are precisely the arbitrary meets of members of $\Omega$: every such meet is fixed by extensivity and monotonicity, and a fixed point has the representation \eqref{eq:kernel}. By \eqref{eq:join}, this family is closed under arbitrary joins as well. Thus $\Fix(\K)$ is the smallest Alexandrov Chang fuzzy topology containing $\Omega$, and $\Omega$ is Alexandrov if and only if $\Omega=\Fix(\K)$.

\section{Point spaces and specialization}\label{sec:recovery}
\subsection{Recovering the fuzzy open family}
Define $\tauM$ to be the topology on $E$ generated by $\{\Pp(U):U\in\Omega\}$. Since
\[
 \Pp(0)=\varnothing,\quad \Pp(1)=E,\quad
 \Pp(U)\cap\Pp(V)=\Pp(U\wedge V),
\]
these sets form a base. Define $\tauQ$ using the sets $\Qq(U)$.

\begin{proposition} \label{prop:Qframe}
The family $\tauQ=\{\Qq(U):U\in\Omega\}$ is already a topology, and
\[
 \Qq:\Omega\longrightarrow\tauQ
\]
is an isomorphism of frames. More precisely, $U$ can be recovered from $\Qq(U)$ by
\begin{equation}\label{eq:q-recovery}
 U(x)=\sup\Bigl(\{t\in[0,1):x_{1-t}\in\Qq(U)\}\cup\{0\}\Bigr).
\end{equation}
\end{proposition}
\begin{proof}
The topology axioms follow from \eqref{eq:Q-ops} and the existence of $0,1\in\Omega$. Formula \eqref{eq:q-recovery} holds because $x_{1-t}\in\Qq(U)$ exactly when $t<U(x)$. Thus $\Qq$ is injective and, by definition of $\tauQ$, surjective. The identities in \eqref{eq:Q-ops} show that it preserves arbitrary joins and finite meets.
\end{proof}
\begin{remark}
The existence of a point representation by an attachment relation and the associated frame map belong to the earlier theory \cite{Guido2010}. Proposition~\ref{prop:Qframe} is included to fix the specific quasi-coincidence formula used below, not as a priority claim.
\end{remark}

The membership topology needs a different recovery argument. For $W\subseteq E$ and $x\in X$, write
\[
 W_x=\{\alpha\in(0,1]:x_\alpha\in W\}.
\]
We view $W_x$ as a subspace of the ordinary interval $(0,1]$.

\begin{lemma}\label{lem:fibers}
If $W\in\tauM$, then each $W_x$ is downward closed. In particular, each vertical section is empty, $(0,a)$, $(0,a]$, or $(0,1]$ for a suitable $a\in(0,1]$.
\end{lemma}
\begin{proof}
Each $\Pp(U)$ has vertical section $(0,U(x)]$ if $U(x)>0$, and the empty section otherwise. Arbitrary unions of such sections are downward closed. Every nonempty downward closed subset of $(0,1]$ has one of the stated forms.
\end{proof}

\begin{theorem} \label{thm:mrecovery}
For $W\subseteq E$, the following are equivalent:
\begin{enumerate}[label=\textup{(\roman*)}]
\item $W=\Pp(U)$ for some $U\in\Omega$;
\item $W$ is open in $\tauM$ and $W_x$ is closed in the usual topology of $(0,1]$ for every $x\in X$.
\end{enumerate}
When these conditions hold, $U$ is unique and is given by
\begin{equation}\label{eq:m-recovery}
 U(x)=\sup(W_x\cup\{0\}).
\end{equation}
Consequently, the membership point topology together with its given level coordinate determines $\Omega$ exactly.
\end{theorem}
\begin{proof}
If $W=\Pp(U)$, then it is open by construction and its vertical sections are closed in $(0,1]$. Conversely, write $W=\bigcup_i\Pp(U_i)$ using the membership base, and set $V=\bigvee_iU_i\in\Omega$. At each $x$, $V(x)=\sup(W_x\cup\{0\})$. By Lemma~\ref{lem:fibers}, $W_x$ is downward closed. If $V(x)>0$, closedness forces its supremum $V(x)$ to belong to $W_x$, so $W_x=(0,V(x)]$; otherwise $W_x$ is empty. Hence $W=\Pp(V)$. The supremum formula also proves uniqueness.
\end{proof}

\begin{corollary}\label{cor:reconstruction}
Let $\Omega$ and $\Lambda$ be Chang fuzzy topologies on the same set $X$, with the same labeled fuzzy-point universe $E$. Then
\[
 \tau_m(\Omega)=\tau_m(\Lambda)\ \Longrightarrow\ \Omega=\Lambda,
\qquad
 \tau_q(\Omega)=\tau_q(\Lambda)\ \Longrightarrow\ \Omega=\Lambda.
\]
\end{corollary}
\begin{proof}
The first implication follows from Theorem~\ref{thm:mrecovery}; the second follows from Proposition~\ref{prop:Qframe}.
\end{proof}

\subsection{Common open sets}\label{sec:common}
For $B\subseteq X$, let $\chi_B$ denote its characteristic function. The following result gives a complete comparison of the two point topologies, without stratification or separation assumptions.

\begin{theorem} \label{thm:intersection}
For every Chang fuzzy topology $\Omega$,
\begin{equation}\label{eq:intersection}
 \tauM\cap\tauQ
 =\{B\times(0,1]: B\subseteq X,\ \chi_B\in\Omega\},
\end{equation}
where $(x,\alpha)$ is identified with the fuzzy point $x_\alpha$.
\end{theorem}
\begin{proof}
Let $W\in\tauM\cap\tauQ$. By Lemma~\ref{lem:fibers}, each $W_x$ is downward closed. Every set $\Qq(U)$ has upward closed vertical sections, and Proposition~\ref{prop:Qframe} shows that all $\tauQ$-open sets have this form. Consequently, each $W_x$ is both downward and upward closed, so it is either empty or all of $(0,1]$. Thus $W=B\times(0,1]$ for some $B\subseteq X$.

Since $W\in\tauQ$, write $W=\Qq(U)$ with $U\in\Omega$. The vertical section of $\Qq(U)$ is empty exactly when $U(x)=0$, and is all of $(0,1]$ exactly when $U(x)=1$. Therefore $U=\chi_B\in\Omega$.

Conversely, if $\chi_B\in\Omega$, then
$\Pp(\chi_B)=\Qq(\chi_B)=B\times(0,1]$, so this set belongs to both topologies.
\end{proof}

\begin{corollary} \label{cor:incomp}
The following conditions are equivalent:
\begin{enumerate}[label=\textup{(\roman*)}]
\item $\tauM=\tauQ$;
\item $\tauM\subseteq\tauQ$;
\item $\tauQ\subseteq\tauM$;
\item every member of $\Omega$ takes only the values $0$ and $1$.
\end{enumerate}
In particular, if $\Omega$ contains a fuzzy open set with a value strictly between $0$ and $1$, then the two point topologies are incomparable.
\end{corollary}
\begin{proof}
Suppose $\tauM\subseteq\tauQ$. For $U\in\Omega$, the set $\Pp(U)$ belongs to both topologies. By Theorem~\ref{thm:intersection}, all its vertical sections are empty or full. This forces $U(x)\in\{0,1\}$ for every $x$.

The same argument starting with $\tauQ\subseteq\tauM$ uses $\Qq(U)$ in place of $\Pp(U)$ and again forces $U$ to be crisp. Conversely, for a crisp open $U$, $\Pp(U)=\Qq(U)$; hence if every open is crisp the families generate the same topology. The remaining implications are immediate.
\end{proof}

\begin{example} \label{ex:onepoint}
Let $X=\{x\}$ and $\Omega=\{0,\underline{1/2},1\}$. Then
\[
 \tauM=\{\varnothing,\ (0,1/2],\ (0,1]\},\quad
 \tauQ=\{\varnothing,\ (1/2,1],\ (0,1]\}.
\]
Here intervals stand for the membership levels of $x$. The intersection is $\{\varnothing,E\}$, and neither topology contains the other. This also shows that the opposite directions of the vertical sections are visible even on a one-point universe.
\end{example}

\subsection{Classical closure and membership specialization}\label{sec:shadow}
The kernel gives a set operator on the fuzzy-point universe:
\begin{equation}\label{eq:Gamma}
 \Gamma(T)=\Pp(\K(J(T)))\qquad(T\subseteq E).
\end{equation}
Besides its order-theoretic closure properties, the operator satisfies a stronger finite-union identity. This gives an exact classical topological interpretation.

\begin{theorem} \label{thm:Gamma}
The operator $\Gamma$ satisfies the Kuratowski closure axioms. Hence it is the closure operator of a classical topology $\tauC$ on $E$. Moreover:
\begin{enumerate}[label=\textup{(\roman*)}]
\item the $\Gamma$-closed sets are precisely $\Pp(C)$ for $C\in\Fix(\K)$;
\item the sets $\Pp(U)$, $U\in\Omega$, form a closed subbase of $\tauC$;
\item for every $T\subseteq E$,
\begin{equation}\label{eq:Gamma-intersect}
 \Gamma(T)=\bigcap\{\Pp(U): U\in\Omega,\ T\subseteq\Pp(U)\}.
\end{equation}
\end{enumerate}
\end{theorem}
\begin{proof}
Extensivity follows from $T\subseteq\Pp(J(T))$ and $J(T)\le\K(J(T))$. Monotonicity is immediate. Since $J(\Pp(C))=C$, idempotence of $\K$ gives
\[
 \Gamma(\Gamma(T))=\Pp\bigl(\K(\K(J(T)))\bigr)=\Gamma(T).
\]
Also $\Gamma(\varnothing)=\Pp(\K(0))=\varnothing$. For $S,T\subseteq E$, we have $J(S\cup T)=J(S)\vee J(T)$, so \eqref{eq:join} and preservation of \emph{binary} joins by $\Pp$ give
\[
 \Gamma(S\cup T)=\Pp(\K(J(S))\vee\K(J(T)))
 =\Gamma(S)\cup\Gamma(T).
\]
Thus $\Gamma$ is a Kuratowski closure operator. If $T$ is fixed then $T=\Pp(C)$ with $C=\K(J(T))$ fixed. Conversely, any $T=\Pp(C)$ with $\K(C)=C$ is fixed. This proves (i).

For (iii), observe that $T\subseteq\Pp(U)$ if and only if $J(T)\le U$. Therefore
\[
 \bigcap_{U\in\Omega,\,J(T)\le U}\Pp(U)
 =\Pp\Bigl(\bigwedge_{U\in\Omega,\,J(T)\le U}U\Bigr)
 =\Pp(\K(J(T)))=\Gamma(T).
\]
In particular, each $\Pp(U)$ is closed and every closed set is an intersection of sets of this form, proving (ii).
\end{proof}

The sets $\Pp(U)$ form an open base for $\tauM$ but a closed subbase for $\tauC$; in general these topologies differ and their specialization orders reverse.

\begin{definition}
For fuzzy points $p=x_\alpha$ and $q=y_\beta$, define
\begin{align*}
 p\preM q&\quad\Longleftrightarrow\quad
 \bigl(p\in\Pp(U)\Rightarrow q\in\Pp(U)\bigr)
 \quad\text{for every }U\in\Omega,\\
 p\preQ q&\quad\Longleftrightarrow\quad
 \bigl(p\in\Qq(U)\Rightarrow q\in\Qq(U)\bigr)
 \quad\text{for every }U\in\Omega.
\end{align*}
These are the specialization preorders of $\tauM$ and $\tauQ$, with the convention $p\le_\tau q$ when every open containing $p$ also contains $q$.
\end{definition}

\begin{proposition} \label{prop:opposite}
For all $p,q\in E$,
\[
 p\le_{\tauC}q\quad\Longleftrightarrow\quad q\preM p.
\]
\end{proposition}
\begin{proof}
By the specialization convention, $p\le_{\tauC}q$ if and only if $p\in\cl_{\tauC}\{q\}$. Theorem~\ref{thm:Gamma} gives $\cl_{\tauC}\{q\}=\Pp(\K(q))$, so this is equivalent to $p\le\K(q)$, or to the assertion that every fuzzy open set containing $q$ also contains $p$.
\end{proof}

\begin{theorem} \label{thm:kernelorder}
For any $A\in\F(X)$ and $y\in X$,
\begin{equation}\label{eq:point-kernel}
 \K(A)(y)=\sup\Bigl(\bigl\{\beta\in(0,1]:\exists x\in X,\ A(x)>0,
 \ x_{A(x)}\preM y_\beta\bigr\}\cup\{0\}\Bigr).
\end{equation}
Consequently, for two Chang fuzzy topologies $\Omega$ and $\Lambda$ on the same set,
\begin{equation}\label{eq:kerneliff}
 \preccurlyeq_m^{\Omega}=\preccurlyeq_m^{\Lambda}
 \quad\Longleftrightarrow\quad \K_\Omega=\K_\Lambda
 \quad\Longleftrightarrow\quad \Fix(\K_\Omega)=\Fix(\K_\Lambda).
\end{equation}
\end{theorem}
\begin{proof}
Every fuzzy set is a pointwise join of the fuzzy points $x_{A(x)}$ for which $A(x)>0$. By \eqref{eq:join},
\[
 \K(A)=\bigvee_{x:A(x)>0}\K(x_{A(x)}).
\]
For $\beta>0$, the condition $x_{A(x)}\preM y_\beta$ is equivalent to $\beta\le\K(x_{A(x)})(y)$ by the definition of $\K$. Taking the supremum over the source points and target levels yields \eqref{eq:point-kernel}.

If the membership preorders agree, \eqref{eq:point-kernel} proves equality of the kernel operators. Conversely, equality of kernels gives equality of $p\preM q\iff q\le\K(p)$. Finally, a closure operator on a complete lattice is uniquely determined by its fixed points, because it sends each element to the least fixed point above it. This gives the last equivalence.
\end{proof}

The next example shows that the preorder need not recover the fuzzy topology

\begin{example} \label{ex:samekernel}
Fix $c\in(0,1)$, let $X=\{x\}$, and identify fuzzy sets with their values at $x$. Define
\[
 \Omega_1=\{0\}\cup(c,1],\qquad
 \Omega_2=\{0\}\cup[c,1].
\]
Both are Chang fuzzy topologies, and $\Omega_1\ne\Omega_2$. Their kernel operators agree:
\[
 \K_i(a)=
 \begin{cases}
 0,&a=0,\\ c,&0<a\le c,\\ a,&a>c.
 \end{cases}
\]
By Theorem~\ref{thm:kernelorder}, their membership specialization preorders agree. Nevertheless, their membership point topologies differ: $\Pp(\underline c)=(0,c]$ is open for $\Omega_2$ but is not open for $\Omega_1$. Their quasi-point topologies also differ by the injectivity in Proposition~\ref{prop:Qframe}.
\end{example}

Thus the membership preorder determines the kernel completion, not necessarily the original Chang topology.

\section{Threshold functions and a representation theorem}\label{sec:profiles}
The threshold profiles below are defined without stratification. Join preservation, composition, membership-order recovery, and kernel reconstruction hold for every Chang fuzzy topology; contraction and diagonal equality require stratification. The matrix classification applies to all Alexandrov Chang topologies. Write $x_t$ for the fuzzy subset with value $t$ at $x$ and zero elsewhere; set $x_0=0$. The open-kernel operator $K=K_\Omega$ was defined in \eqref{eq:kernel}.

\begin{definition}\label{def:h}
For $x,y\in X$ define $h_{xy}:[0,1]\to[0,1]$ by
\begin{equation}\label{eq:defh}
 h_{xy}(t)=K(x_t)(y)=\inf\{U(y):U\in\Omega,\ U(x)\ge t\}.
\end{equation}
\end{definition}
The set over which the infimum is taken is nonempty: the constant $1$ is open. Although this formula resembles the specialization construction of established fuzzy-order theory \cite{LaiZhang2006,Fang2007}, keeping the complete function $t\mapsto h_{xy}(t)$, rather than one number, lets us compare different membership grades without imposing a chosen implication.

\begin{theorem} \label{thm:matrixproperties}
For every Chang fuzzy topology the functions $h_{xy}$ satisfy the following properties, with stratification required only for the stronger identities in \textup{(ii)}:
\begin{enumerate}[label=\textup{(\roman*)}]
\item $h_{xy}(0)=0$ and $h_{xy}(\sup_{i\in I} t_i)=\sup_{i\in I}h_{xy}(t_i)$ for every family, including the empty family;
\item $h_{xx}(t)\ge t$ for every $x,t$. If $\Omega$ is stratified, then $h_{xx}(t)=t$ and $0\le h_{xy}(t)\le t$ for every $x,y,t$;
\item for all $x,y,z\in X$ and $t\in[0,1]$,
\begin{equation}\label{eq:comp}
 h_{yz}(h_{xy}(t))\le h_{xz}(t);
\end{equation}
\item the membership point order is given exactly by
\begin{equation}\label{eq:ordermatrix}
 x_\alpha\preM y_\beta\quad\Longleftrightarrow\quad \beta\le h_{xy}(\alpha);
\end{equation}
\item for every $A\in[0,1]^X$ and $y\in X$,
\begin{equation}\label{eq:matrixkernel}
 K(A)(y)=\sup_{x\in X}h_{xy}(A(x)).
\end{equation}
\end{enumerate}
\end{theorem}
\begin{proof}
The kernel preserves arbitrary joins. Since $x_{\sup_i t_i}=\bigvee_i x_{t_i}$, evaluating at $y$ proves (i). Extensivity gives $K(x_t)(x)\ge t$ without stratification. If $\Omega$ is stratified, the constant fuzzy open $\underline t$ contains $x_t$, giving $K(x_t)\le\underline t$ and hence the stronger identities in (ii). For (iii), the fuzzy subset $y_{h_{xy}(t)}$ (zero if its grade is zero) lies below $K(x_t)$. Applying the monotone operator $K$ and idempotence gives $K(y_{h_{xy}(t)})\le K(x_t)$. Evaluating at $z$ yields \eqref{eq:comp}. Assertion (iv) follows directly from the definition of the membership specialization order. Finally, $A=\bigvee_{x\in X}x_{A(x)}$. Applying the arbitrary-join identity and evaluating at $y$ gives (v).
\end{proof}

The preceding conditions characterize admissible threshold matrices.
This is a coordinate instance of established join-preserving
closure-operator and quantale-enriched relation theory
\cite{LaiZhang2006,Fang2007,Kim2014join,OhKim2021,JarvinenKortelainen2007},
not a new abstract correspondence. The matrix representation retains
unequal-grade comparisons without assuming scalar compatibility,
as required for the realization-fiber results.

\begin{theorem} \label{thm:matrixrep}
Let $X\ne\varnothing$. Suppose maps $h_{xy}:[0,1]\to[0,1]$, for $x,y\in X$, satisfy conditions \textup{(i)} and \textup{(iii)} of Theorem~\ref{thm:matrixproperties}, together with $h_{xx}(t)\ge t$ for every $x,t$. Define
\begin{equation}\label{eq:KH}
 K_h(A)(y)=\sup_{x\in X}h_{xy}(A(x)).
\end{equation}
Then $K_h$ is an extensive, monotone, idempotent, arbitrary-join-preserving operator on $[0,1]^X$. Its fixed-point family
\[
 \Omega_h=\{A\in[0,1]^X:K_h(A)=A\}
\]
is an Alexandrov Chang fuzzy topology. Its open kernel is $K_h$, and $K_h(x_t)(y)=h_{xy}(t)$. Conversely, every Alexandrov Chang fuzzy topology occurs this way and determines the matrix uniquely.

Moreover, $\Omega_h$ is stratified if and only if $h_{xy}(t)\le t$ for every $x,y,t$. Under this additional condition, $h_{xx}(t)=t$ for every $x,t$.
\end{theorem}
\begin{proof}
Each $h_{xy}$ is monotone because it preserves arbitrary suprema. Hence $K_h$ is monotone and preserves arbitrary joins by interchanging the two suprema. The diagonal inequality gives $A(y)\le h_{yy}(A(y))\le K_h(A)(y)$, so $K_h$ is extensive. For idempotence,
\begin{align*}
 K_h^2(A)(z)
 &=\sup_{y\in X}h_{yz}\left(\sup_{x\in X}h_{xy}(A(x))\right)\\
 &=\sup_{x,y\in X}h_{yz}(h_{xy}(A(x)))
 \le \sup_xh_{xz}(A(x))=K_h(A)(z).
\end{align*}
Extensivity gives the reverse inequality. The fixed points of a closure operator are closed under arbitrary meets; join preservation also gives closure under arbitrary joins, including the empty join. Thus $\Omega_h$ is an Alexandrov Chang topology. Since $K_h(A)$ is the least fixed point above $A$, its open kernel is $K_h$. Evaluating at $x_t$ gives $h_{xy}(t)$, because $h_{zy}(0)=0$. The converse and uniqueness follow from Theorem~\ref{thm:matrixproperties} and the equality of an Alexandrov topology with its kernel completion.

Finally, $K_h(\underline t)(y)=\sup_x h_{xy}(t)$. Every constant is fixed if and only if all entries are bounded above by the identity: the forward implication follows from the supremum formula, and the reverse from the diagonal inequality. This also forces diagonal equality.
\end{proof}

\begin{remark}[Operator-matrix perspective]
Let $\operatorname{Sup}([0,1])$ denote the arbitrary-join-preserving self-maps of $[0,1]$, ordered pointwise. Pointwise suprema and composition make this a unital (generally noncommutative) quantale; composition distributes over arbitrary suprema because the maps preserve them. The profiles $h_{xy}$ may thus be viewed as matrix entries in this operator quantale. Condition \eqref{eq:comp} is the corresponding transitivity inequality, and the diagonal entries dominate the identity (and equal it in the stratified case). This observation situates Theorem~\ref{thm:matrixrep} within general operator-enriched representation theory; the paper's later questions concern the information lost by replacing each operator entry by a single G\"odel grade.
\end{remark}

\begin{corollary}[Matrix form of membership recovery]\label{cor:orderrec}
For two Chang fuzzy topologies on the same labeled fuzzy-point universe,
equality of their membership-point preorders is equivalent to equality
of their threshold matrices. Consequently, it is also equivalent to
equality of their open-kernel operators and Alexandrov completions.
If both topologies are Alexandrov, they are equal.
\end{corollary}
\begin{proof}
Equation~\eqref{eq:ordermatrix} recovers each preorder from its matrix;
conversely, each matrix entry is the supremum of the target grades
permitted by that preorder. Equation~\eqref{eq:matrixkernel}
reconstructs the kernel, whose fixed-point family is the Alexandrov
completion. The last assertion follows immediately.
\end{proof}

\section{Exactly what G\"odel specialization forgets}\label{sec:loss}
Let $\to_G$ denote the G\"odel implication on $[0,1]$:
\[
 u\to_Gv=\begin{cases}1,&u\le v,\\v,&u>v.\end{cases}
\]
The graded specialization relation is
\begin{equation}\label{eq:rdef}
 r_\Omega(x,y)=\inf_{U\in\Omega}\bigl(U(x)\to_GU(y)\bigr).
\end{equation}
Reflexivity is immediate. For every $U\in\Omega$, residuation gives $\min\{U(x),r_\Omega(x,y)\}\le U(y)$, and hence
\[
 \min\{U(x),r_\Omega(x,y),r_\Omega(y,z)\}\le U(z).
\]
Applying residuation again and taking the infimum over $U$ proves minimum-transitivity. The relation and its fuzzy-upper-set representation are standard \cite{LaiZhang2006,Fang2007}. The next result recovers it from the full membership order without any stratification assumption.

\begin{theorem} \label{thm:gradedprofile}
For every Chang fuzzy topology and all $x,y\in X$,
\begin{equation}\label{eq:gradedprofile}
 r_\Omega(x,y)=\inf_{t\in[0,1]}\bigl(t\to_G h_{xy}(t)\bigr).
\end{equation}
In particular, the membership-point specialization preorder determines the G\"odel relation, without stratification.
\end{theorem}
\begin{proof}
Put $R=r_\Omega(x,y)$. For every $U\in\Omega$ with $U(x)\ge t$, residuation gives $\min\{t,R\}\le\min\{U(x),R\}\le U(y)$. Taking the infimum over these opens yields $\min\{t,R\}\le h_{xy}(t)$, so $R\le t\to_G h_{xy}(t)$ for every $t$. Conversely, let $S$ be the right side of \eqref{eq:gradedprofile}. For any $U\in\Omega$, taking $t=U(x)$ gives $h_{xy}(t)\le U(y)$. Monotonicity in the second argument of G\"odel implication yields
\[
 S\le U(x)\to_G h_{xy}(U(x))\le U(x)\to_G U(y).
\]
Taking the infimum over $U$ proves $S\le R$. The last assertion follows from Corollary~\ref{cor:orderrec}.
\end{proof}

\begin{theorem} \label{thm:firstdeparture}
Let $\Omega$ be stratified. For every $x,y\in X$,
\begin{equation}\label{eq:firstdeparture}
 r_\Omega(x,y)=\inf\bigl(\{t\in(0,1]:h_{xy}(t)<t\}\cup\{1\}\bigr).
\end{equation}
Equivalently, for $c\in[0,1]$,
\begin{equation}\label{eq:equalall}
 c\le r_\Omega(x,y)
 \quad\Longleftrightarrow\quad
 x_t\preM y_t\text{ for every }0<t\le c.
\end{equation}
Equivalently,
\begin{equation}\label{eq:initialendpoint}
 r_\Omega(x,y)=\sup\{c\in[0,1]:h_{xy}(t)=t\text{ for all }0\le t\le c\}.
\end{equation}
Thus the graded value is the endpoint of the initial interval of identity agreement. The comparison holds at this endpoint (when positive); it need not fail there, and it may hold again at larger levels.
\end{theorem}
\begin{proof}
Let $d$ denote the right side of \eqref{eq:firstdeparture}. If $U(x)>U(y)=b$, then for every $t$ with $b<t\le U(x)$, the open set $U$ occurs in the infimum defining $h_{xy}(t)$, and $h_{xy}(t)\le b<t$. Hence $d\le b$. Taking the infimum over all such $U$ gives $d\le r_\Omega(x,y)$ (if there is no violating open, both equal $1$).

Conversely, suppose $h_{xy}(t)<t$. Since $h_{xy}(t)$ is the infimum of the numbers $U(y)$ for opens with $U(x)\ge t$, there exists such a $U$ with $U(y)<t$. Then $U(x)>U(y)$, so $r_\Omega(x,y)\le U(y)<t$. Taking the infimum over all such $t$ yields $r_\Omega(x,y)\le d$. This proves \eqref{eq:firstdeparture}.

If $c>r_\Omega(x,y)$, the infimum formula guarantees a departure at some $t<c$, so the right-hand side of \eqref{eq:equalall} fails. If $c<r_\Omega(x,y)$ there is no departure at any $t\le c$, giving \eqref{eq:equalall} away from the endpoint. At $c=r_\Omega(x,y)>0$, suppose $h_{xy}(c)<c$. For every $t<c$ no departure is possible, so $h_{xy}(t)=t$. Since $h_{xy}$ preserves arbitrary suprema and $c=\sup_{t<c}t$, this implies $h_{xy}(c)=c$, a contradiction. The case $c=0$ is vacuous. Combining with \eqref{eq:ordermatrix} proves \eqref{eq:equalall}.
\end{proof}

\begin{remark}
The equivalence \eqref{eq:equalall} is compatible with the standard threshold interpretation of G\"odel implication. The added information in our comparison is the full function $h_{xy}$ and the exact description of the many matrices that can have the same first-departure levels. No novelty is claimed for the basic graded specialization relation itself.
\end{remark}

\section{Scalar compatibility and recovery from graded specialization}\label{sec:scalar}
The representation in Theorem~\ref{thm:matrixrep} assumes arbitrary-meet closure but \emph{not} compatibility with minimum by constants or with scalar G\"odel implication. These extra algebraic conditions occur in fuzzy-preorder representations \cite{LaiZhang2006,Fang2007} and approximation-operator treatments \cite{OhKim2021}, but they do not follow from our Alexandrov definition. The next results isolate this difference using threshold functions; imposing the stronger convention from the outset would incorrectly collapse the non-scalar fibers of Section~\ref{sec:twopoint}.

\begin{lemma} \label{lem:thresholdbound}
For any stratified Chang fuzzy topology, write $r=r_\Omega$ and $h=h^\Omega$. Then for all $x,y\in X$ and $t\in[0,1]$,
\begin{equation}\label{eq:lowerbound}
 \min\{t,r(x,y)\}\le h_{xy}(t)\le t.
\end{equation}
In particular, $h_{xy}(t)=t$ whenever $t\le r(x,y)$, and $r(x,y)\le h_{xy}(1)$.
\end{lemma}
\begin{proof}
The upper bound follows from stratification. If $r(x,y)=0$, the lower bound is immediate. Suppose $r(x,y)>0$. By Theorem~\ref{thm:firstdeparture}, $h_{xy}(s)=s$ for every $s<r(x,y)$. Because $h_{xy}$ preserves arbitrary suprema, $h_{xy}(r(x,y))=r(x,y)$, including when $r(x,y)=1$. By monotonicity, for $t\ge r(x,y)$ we have $h_{xy}(t)\ge r(x,y)$. This gives \eqref{eq:lowerbound}.
\end{proof}

\begin{theorem} \label{thm:scalar}
Let $\Omega$ be a stratified Chang fuzzy topology with kernel $K$, threshold functions $h_{xy}$ and G\"odel specialization relation $r=r_\Omega$. The following conditions are equivalent:
\begin{enumerate}[label=\textup{(\roman*)}]
\item $K(\underline c\wedge A)=\underline c\wedge K(A)$ for all $c\in[0,1]$ and $A\in[0,1]^X$;
\item $h_{xy}(t)=\min\{t,h_{xy}(1)\}$ for all $x,y\in X$ and $t\in[0,1]$;
\item $h_{xy}(1)=r(x,y)$ for all $x,y\in X$;
\item $h_{xy}(t)=\min\{t,r(x,y)\}$ for all $x,y\in X$ and $t\in[0,1]$.
\end{enumerate}
When these conditions hold, the \emph{kernel completion} is determined by the G\"odel relation:
\begin{equation}\label{eq:gradedkernel}
 K(A)(y)=\sup_{x\in X}\min\{A(x),r(x,y)\}.
\end{equation}
The equality $\Omega=\Fix(K)$, and hence reconstruction of the \emph{original} topology from $r$, additionally requires the Alexandrov property.
\end{theorem}
\begin{proof}
For (i)$\Rightarrow$(ii), use $x_t=\underline t\wedge x_1$ and evaluate the identity at $y$. If (ii) holds and $s=h_{xy}(1)$, the first-departure formula says that $r(x,y)=s$: the function $t\mapsto\min\{t,s\}$ first falls below $t$ immediately above $s$ when $s<1$, and never does when $s=1$. Hence (ii) implies (iii).

Suppose (iii). By Lemma~\ref{lem:thresholdbound}, $h_{xy}(t)=t$ for $t\le r(x,y)$. For $t>r(x,y)$, monotonicity gives $h_{xy}(t)\ge r(x,y)$, while $h_{xy}(t)\le h_{xy}(1)=r(x,y)$. Thus (iv) holds. Finally, using (iv) and the matrix formula \eqref{eq:matrixkernel},
\begin{align*}
 K(\underline c\wedge A)(y)
 &=\sup_x\min\{c,A(x),r(x,y)\}\\
 &=\min\left\{c,\sup_x\min\{A(x),r(x,y)\}\right\}\\
 &=\min\{c,K(A)(y)\},
\end{align*}
where the second identity is distributivity of finite meets over arbitrary joins in $[0,1]$. This gives (i). The formula \eqref{eq:gradedkernel} follows from (iv) and \eqref{eq:matrixkernel}.
\end{proof}

\begin{proposition} \label{prop:residualboundary}
Suppose that $\Omega$ is a stratified Alexandrov Chang fuzzy topology. Then the equivalent conditions of Theorem~\ref{thm:scalar} hold if and only if
\begin{equation}\label{eq:residualclosure}
 (\underline c\to_G U)\in\Omega
 \qquad(c\in[0,1],\ U\in\Omega),
\end{equation}
where $\underline c\to_G U$ is formed pointwise. Thus arbitrary-meet closure and containment of all constants do not entail the residual-stability condition used in standard fuzzy-preorder--Alexandrov correspondences.
\end{proposition}
\begin{proof}
Assume first that $K(\underline c\wedge A)=\underline c\wedge K(A)$ for all $c,A$, and fix $U\in\Omega$. Let $B=\underline c\to_G U$. By residuation, $\underline c\wedge B\le U$. Thus
\[
 \underline c\wedge K(B)=K(\underline c\wedge B)\le K(U)=U.
\]
Residuation gives $K(B)\le\underline c\to_G U=B$. By extensivity $K(B)=B$; since $\Omega$ is Alexandrov, $B\in\Omega$.

Conversely, assume \eqref{eq:residualclosure}. Since $K$ is monotone and the constant $\underline c$ is open,
\[
 K(\underline c\wedge A)\le K(A)\wedge\underline c.
\]
If $U\in\Omega$ and $\underline c\wedge A\le U$, then $A\le\underline c\to_G U\in\Omega$, so
$K(A)\le\underline c\to_G U$ and hence $\underline c\wedge K(A)\le U$.
Taking the infimum over all such $U$ gives
$\underline c\wedge K(A)\le K(\underline c\wedge A)$, as required.
\end{proof}

\begin{remark} \label{rem:scalarprior}
Lai--Zhang~\cite[Theorems~3.15 and~4.8]{LaiZhang2006} characterize fuzzy upper-set topologies using, among other conditions, stability under scalar implication. For the minimum t-norm, \eqref{eq:gradedkernel} is the familiar relational upper-set closure, not a new representation formula. Proposition~\ref{prop:residualboundary} and Theorem~\ref{thm:scalar} express this extra requirement directly in the threshold notation. In contrast, our Alexandrov class requires arbitrary joins and meets and all constants, but not scalar residual stability.

Kim~\cite[Theorem~2.9]{Kim2014join} studies operators preserving both arbitrary joins and scalar multiplication in a setting with strong negation, so the results there do not specialize directly to the non-involutive G\"odel negation. The analogous scalar law for minimum is condition~(i) of Theorem~\ref{thm:scalar}. It is strictly stronger than preservation of joins alone: for the family in Example~\ref{thm:continuum}, with $b<a<1$ and $b<t\le a$, one has $K(x_t)(y)=b<t=\min\{t,K(x_1)(y)\}$. Related distinctions between meet closure and scalar compatibility occur in fuzzy rough approximation theory \cite[Remark~2 and Theorem~26]{XuZouLi2023}.
\end{remark}

\begin{theorem} \label{thm:gradedmax}
Let $r:X\times X\to[0,1]$ be a reflexive, minimum-transitive fuzzy preorder, that is,
\[
 r(x,x)=1,\qquad \min\{r(x,y),r(y,z)\}\le r(x,z).
\]
Put
\begin{equation}\label{eq:Omegar}
 \Omega_r=\bigl\{A\in[0,1]^X:\min\{A(x),r(x,y)\}\le A(y)\ \text{for every }x,y\in X\bigr\}.
\end{equation}
Then $\Omega_r$ is a stratified Alexandrov Chang fuzzy topology with G\"odel specialization $r$. More strongly, \emph{every Chang fuzzy topology} $\Omega$ (even if neither stratified nor Alexandrov) satisfying $r_\Omega=r$ is contained in $\Omega_r$. For a \emph{stratified Alexandrov} $\Omega$ realizing $r$, the following are equivalent:
\begin{enumerate}[label=\textup{(\roman*)}]
\item $\Omega=\Omega_r$;
\item its kernel commutes with meets with constants;
\item its threshold functions satisfy $h_{xy}(1)=r(x,y)$ for all $x,y$.
\end{enumerate}
\end{theorem}
\begin{proof}
Put $h^r_{xy}(t)=\min\{t,r(x,y)\}$. Each $h^r_{xy}$ preserves arbitrary suprema and is bounded above by the identity. Reflexivity gives $h^r_{xx}=\operatorname{id}$; minimum-transitivity gives
\[
 h^r_{yz}(h^r_{xy}(t))=\min\{t,r(x,y),r(y,z)\}\le h^r_{xz}(t).
\]
Theorem~\ref{thm:matrixrep} constructs a stratified Alexandrov topology with kernel
\[
 K_r(A)(y)=\sup_{x\in X}\min\{A(x),r(x,y)\}.
\]
Its fixed points are exactly \eqref{eq:Omegar}, and Theorem~\ref{thm:firstdeparture} gives its graded specialization $r$. This is the standard fuzzy-upper-set construction \cite{LaiZhang2006,Fang2007}.

Let $\Omega$ be any Chang fuzzy topology with $r_\Omega=r$, and let $U\in\Omega$. By definition of graded specialization,
\[
 r(x,y)\le U(x)\to_G U(y)
 \quad\Longrightarrow\quad
 \min\{U(x),r(x,y)\}\le U(y)
\]
for all $x,y\in X$. Hence $U\in\Omega_r$, proving $\Omega\subseteq\Omega_r$ \emph{without} stratification or the Alexandrov property. Now restrict to stratified Alexandrov $\Omega$. If $\Omega=\Omega_r$, its kernel has the form \eqref{eq:gradedkernel} and commutes with meets with constants. Conversely, if the kernel commutes with these meets, Theorem~\ref{thm:scalar} gives $K_\Omega=K_r$, and both topologies equal their fixed-point families because they are Alexandrov. The equivalence with (iii) is Theorem~\ref{thm:scalar}.
\end{proof}

\begin{remark}
The existence of the fuzzy-upper-set topology and its connection with fuzzy specialization are established theory \cite{LaiZhang2006,Fang2007}. The inclusion above is presented to identify precisely where the larger family of ordinary Chang Alexandrov topologies fits relative to that classical construction, rather than to assert a new fuzzy-order representation theorem.
\end{remark}

\section{Realizations of a prescribed graded relation}\label{sec:generalfibre}
For a reflexive, minimum-transitive relation $r:X\times X\to[0,1]$, let $\mathfrak T_r$ be the set of all stratified Chang fuzzy topologies on the fixed nonempty set $X$ with G\"odel specialization $r$, ordered by inclusion. Write $\mathfrak A_r\subseteq\mathfrak T_r$ for its Alexandrov members. By Theorem~\ref{thm:gradedmax}, both families contain $\Omega_r$ as their greatest element. We now study both families by a single construction on the original open sets, without assuming arbitrary-meet closure.

\begin{lemma} \label{lem:enlarge}
Let $\Omega$ be a stratified Chang fuzzy topology and $0\le\lambda<1$. Define
\begin{equation}\label{eq:uniformrestriction}
 \Omega^\lambda=\{U\in\Omega:U(y)\ge\lambda U(x)\text{ for all }x,y\in X\}.
\end{equation}
Then $\Omega^\lambda$ is stratified, and is Alexandrov whenever $\Omega$ is Alexandrov. Its kernel is given by
\begin{equation}\label{eq:uniformkernel}
 K_{\Omega^\lambda}(A)(y)=\max\{K_\Omega(A)(y),\lambda\sup_{x\in X}A(x)\}.
\end{equation}
In particular, its threshold profiles satisfy
\begin{equation}\label{eq:enlarge}
 h^\lambda_{xy}(t)=\max\{h_{xy}(t),\lambda t\},
\end{equation}
and $r_{\Omega^\lambda}=r_\Omega$. If $\lambda_1\le\lambda_2$, then $\Omega^{\lambda_2}\subseteq\Omega^{\lambda_1}\subseteq\Omega$. These inclusions are strict whenever the corresponding kernels differ.
\end{lemma}
\begin{proof}
The defining inequalities hold for every constant and are preserved by arbitrary joins and finite meets: multiplication by $\lambda$ commutes with both operations on $[0,1]$ (for the empty meet, the inequality holds directly for $1$). For the same reason, arbitrary nonempty meets preserve the inequalities, so $\Omega^\lambda$ is Alexandrov if $\Omega$ is. The inclusion assertions follow from monotonicity in $\lambda$.

Fix $A\in[0,1]^X$ and let $a=\sup_x A(x)$. Since $\Omega^\lambda\subseteq\Omega$, its kernel dominates $K_\Omega$. Each $V\in\Omega^\lambda$ with $A\le V$ also satisfies $V(y)\ge\lambda V(x)\ge\lambda A(x)$ for all $x$, hence $V(y)\ge\lambda a$. Taking the infimum of these opens gives the lower bound in \eqref{eq:uniformkernel}.

For the reverse bound, take any $U\in\Omega$ with $A\le U$. Stratification and the topology axioms imply
\[
 V=(U\wedge\underline a)\vee\underline{\lambda a}\in\Omega.
\]
Since $\lambda a\le V(x)\le a$ for every $x$, one has $V\in\Omega^\lambda$. Also $A\le V$, and $V(y)\le\max\{U(y),\lambda a\}$. Consequently
\[
 K_{\Omega^\lambda}(A)(y)\le\inf_{\substack{U\in\Omega\\A\le U}}\max\{U(y),\lambda a\}
 =\max\{K_\Omega(A)(y),\lambda a\},
\]
where the last equality uses distributivity of maximum with a fixed scalar over nonempty infima in $[0,1]$. This proves \eqref{eq:uniformkernel}; taking $A=x_t$ gives \eqref{eq:enlarge}. For $t>0$, $\lambda t<t$, so $h^\lambda_{xy}(t)<t$ if and only if $h_{xy}(t)<t$. Theorem~\ref{thm:firstdeparture} gives equality of the graded relations. Equal open families have equal kernels, which proves the strictness assertion.
\end{proof}

\begin{remark}
For an admissible stratified matrix $h$ of Theorem~\ref{thm:matrixrep}, applying Lemma~\ref{lem:enlarge} to $\Omega_h$ shows that the matrix $h^\lambda$ in \eqref{eq:enlarge} is again admissible and that $\Omega_{h^\lambda}=(\Omega_h)^\lambda$. Thus the construction restricts the Alexandrov family as well as the full stratified family.
\end{remark}

\begin{theorem} \label{thm:generalfibre}
Let $X$ be nonempty and $r$ a reflexive, minimum-transitive fuzzy preorder on $X$. Then:
\begin{enumerate}[label=\textup{(\roman*)}]
\item Both $\mathfrak T_r$ and $\mathfrak A_r$ have greatest element $\Omega_r$ given by \eqref{eq:Omegar}.
\item If $r\equiv1$, both families consist of the single topology of all constant fuzzy sets.
\item If $r\not\equiv1$, for every $\Omega\in\mathfrak T_r$ there is a family $(\Omega^\lambda)_{\lambda\in I}$, indexed by a nonempty open interval $I\subseteq[0,1)$, such that
\[
 \Omega^{\lambda_2}\subsetneq\Omega^{\lambda_1}\subsetneq\Omega
 \qquad(\lambda_1<\lambda_2,\ \lambda_1,\lambda_2\in I),
\]
and every $\Omega^\lambda$ belongs to $\mathfrak T_r$. If $\Omega\in\mathfrak A_r$, all these topologies belong to $\mathfrak A_r$. Thus neither family has an inclusion-minimal member, and strictly below every member there is a chain of cardinality $2^{\aleph_0}$ within its respective family. The interval $I$ can be chosen with supremum $1$, in which case
\begin{equation}\label{eq:chainintersection}
 \bigcap_{\lambda\in I}\Omega^\lambda=\{\underline t:t\in[0,1]\}.
\end{equation}
The constructed chain has no lower bound in $\mathfrak T_r$. If $\Omega\in\mathfrak A_r$, the chain lies entirely in $\mathfrak A_r$ and has no lower bound there either.
\item If $X$ is finite or countably infinite and $r\not\equiv1$, then $|\mathfrak A_r|=2^{\aleph_0}$.
\end{enumerate}
\end{theorem}
\begin{proof}
Part (i) follows from Theorem~\ref{thm:gradedmax}. If $r\equiv1$, that theorem implies every realizing topology is contained in $\Omega_r$, which in this case consists of constants only. Stratification forces equality, proving (ii).

For (iii), let $\Omega\in\mathfrak T_r$ and let $h$ be its threshold matrix. Since $r\not\equiv1$, Theorem~\ref{thm:firstdeparture} gives $x,y\in X$ and $t>0$ with $h_{xy}(t)<t$. Put $I=(h_{xy}(t)/t,1)$. For each $\lambda\in I$, Lemma~\ref{lem:enlarge} gives a realization $\Omega^\lambda\subseteq\Omega$ with $h^\lambda_{xy}(t)=\lambda t>h_{xy}(t)$. Thus each inclusion in $\Omega$ is strict. The values $\lambda t$ are strictly increasing in $\lambda$, so the kernels are distinct; the inclusion assertion of the lemma gives strict reverse inclusion of the open families. This constructs the required continuum-sized chain strictly below $\Omega$. The preservation of the Alexandrov property in the lemma proves the corresponding assertion for $\mathfrak A_r$. For this choice of $I$, if $U$ lies in the intersection, then $U(y)\ge\lambda U(x)$ for every $\lambda\in I$ and all $x,y$. Taking the supremum over $\lambda$ gives $U(y)\ge U(x)$; interchanging $x,y$ yields equality. Conversely, all constants belong to every restricted topology. This proves \eqref{eq:chainintersection}. A common lower bound within $\mathfrak T_r$ would necessarily equal the constant topology by stratification, but its specialization is identically $1$, contradicting $r\not\equiv1$. The same reasoning applies to $\mathfrak A_r$.

For (iv), a join-preserving map $f:[0,1]\to[0,1]$ satisfies $f(0)=0$ and
\[
 f(t)=\sup\{f(q):q\in\mathbb Q,\ 0\le q<t\}\qquad(t>0).
\]
There are at most $2^{\aleph_0}$ assignments of values to countably many rational inputs. If $X$ is at most countable, an admissible matrix is determined by at most countably many such maps, so there are at most $(2^{\aleph_0})^{\aleph_0}=2^{\aleph_0}$ matrices. Theorem~\ref{thm:matrixrep} identifies the Alexandrov topologies with their matrices. The chain in (iii) supplies the matching lower bound.
\end{proof}

\begin{remark}
Starting at the greatest element $\Omega_r$, the construction yields the explicit matrices
\begin{equation}\label{eq:uniformfamily}
 m^{r,\lambda}_{xy}(t)=\max\{\min(t,r(x,y)),\lambda t\}.
\end{equation}
If $r(x_0,y_0)=d<1$, the entries $m^{r,\lambda}_{x_0y_0}(1)=\lambda$ for $d<\lambda<1$ directly distinguish a continuum-sized chain of realizations. The absence of a least realization was already observed in particular examples by Lai--Zhang~\cite[Example~4.24]{LaiZhang2006}. Theorem~\ref{thm:generalfibre} applies to every nontrivial prescribed relation and constructs such a chain strictly below each stratified realization. Its cardinality upper bound is specifically for Alexandrov realizations: counting matrices alone does not bound the number of distinct open families with one kernel. The descending chain has no lower bound in the fiber; thus neither nontrivial fiber is filtered-complete (that is, admits an infimum for every nonempty downward directed subset). Theorem~\ref{thm:finiteallcard} below obtains the full finite-set cardinality by a different argument, using quasi-specialization.
\end{remark}

\section{The two-point realization fiber}\label{sec:twopoint}
The general theorem does not identify the entire inclusion order of a fiber. On two points the matrix conditions simplify sufficiently to give an explicit product classification.

\begin{theorem} \label{thm:twopointclass}
Let $X=\{x,y\}$. Stratified Alexandrov Chang fuzzy topologies on $X$ are in one-to-one correspondence with ordered pairs $(p,q)$ of maps $[0,1]\to[0,1]$ such that each preserves arbitrary suprema and satisfies
\[
 0\le p(t)\le t,\qquad 0\le q(t)\le t \quad(t\in[0,1]).
\]
The correspondence sends a topology to
\[
 p(t)=K(x_t)(y),\qquad q(t)=K(y_t)(x),
\]
and its inverse consists of all pairs $(u,v)\in[0,1]^2$ satisfying
\begin{equation}\label{eq:fixedpq}
 u\ge q(v),\qquad v\ge p(u).
\end{equation}
For this topology the G\"odel specialization values are exactly
\begin{equation}\label{eq:rpq}
 r(x,y)=\inf(\{t>0:p(t)<t\}\cup\{1\}),\qquad
 r(y,x)=\inf(\{t>0:q(t)<t\}\cup\{1\}).
\end{equation}
\end{theorem}
\begin{proof}
For a stratified Alexandrov topology, Theorem~\ref{thm:matrixproperties} gives the required properties of $p,q$. Conversely, take arbitrary $p,q$ with the stated properties and define $h_{xx}=h_{yy}=\operatorname{id}$, $h_{xy}=p$, $h_{yx}=q$. The composition inequalities of Theorem~\ref{thm:matrixproperties} involving three points reduce either to $p\le\operatorname{id}$, $q\le\operatorname{id}$, $q\circ p\le\operatorname{id}$, or $p\circ q\le\operatorname{id}$. All follow from $p(t),q(t)\le t$ and monotonicity. The matrix representation theorem therefore supplies a unique stratified Alexandrov topology. Its kernel is
\[
 K(u,v)=\bigl(\max\{u,q(v)\},\max\{v,p(u)\}\bigr),
\]
so its fixed points are exactly \eqref{eq:fixedpq}. These formulas also show injectivity of the correspondence. The expressions in \eqref{eq:rpq} follow from Theorem~\ref{thm:firstdeparture}.
\end{proof}

\begin{corollary} \label{cor:twoindependent}
Every pair $(b,c)\in[0,1]^2$ occurs as $(r(x,y),r(y,x))$ for a stratified Alexandrov Chang topology on two points. For a fixed pair of graded values, the corresponding topologies are described exactly by those pairs of join-preserving contractions $(p,q)$ satisfying the two first-departure conditions in \eqref{eq:rpq}.
\end{corollary}
\begin{proof}
Take $p(t)=\min\{t,b\}$ and $q(t)=\min\{t,c\}$. Both preserve arbitrary suprema and are bounded above by the identity. Apply Theorem~\ref{thm:twopointclass}; their first departures are $b$ and $c$, with the empty departure set interpreted as value $1$.
\end{proof}

\begin{example} \label{thm:continuum}
Fix $b\in(0,1)$ and let $X=\{x,y\}$. For each $a\in(b,1)$ define
\begin{equation}\label{eq:Omegaab}
 \Omega_{a,b}=\{\underline t:t\in[0,1]\}\;\cup\;
 \{U_{u,v}: b\le v<u\le a\},
\end{equation}
where $\underline t=(t,t)$ and $U_{u,v}(x)=u$, $U_{u,v}(y)=v$. Then:
\begin{enumerate}[label=\textup{(\roman*)}]
\item each $\Omega_{a,b}$ is a stratified Alexandrov Chang fuzzy topology;
\item its threshold functions are
\begin{equation}\label{eq:hab}
 h_{xy}^{a,b}(t)=\begin{cases}
 t,&0\le t\le b,\\b,&b<t\le a,\\t,&a<t\le1,
 \end{cases}
 \qquad h_{yx}^{a,b}(t)=t,
\end{equation}
with diagonal functions equal to the identity;
\item its G\"odel specialization is independent of $a$:
\[
 r_{a,b}(x,y)=b,\qquad r_{a,b}(y,x)=1;
\]
\item if $b<a_1<a_2<1$, then $\Omega_{a_1,b}\subsetneq\Omega_{a_2,b}$, their kernel operators differ, and their unequal-grade membership preorders differ.
\end{enumerate}
In particular, a single fixed graded specialization relation is shared by continuum many distinct stratified Alexandrov fuzzy topologies and continuum many distinct full fuzzy-point orders, even on a two-point set.
\par\medskip\noindent\textit{Verification.}
Let $p:[0,1]\to[0,1]$ be the piecewise function in \eqref{eq:hab} and put $q(t)=t$. The function $p$ is monotone and satisfies $p(0)=0$ and $p(t)=\sup_{s<t}p(s)$ for every $t>0$, including $t=a$ and $t=b$. A monotone function with these properties preserves arbitrary suprema on $[0,1]$. Both $p$ and $q$ are bounded above by the identity. By Theorem~\ref{thm:twopointclass}, there is a unique stratified Alexandrov topology with threshold matrix $(p,q)$; its open fuzzy sets are exactly the pairs $(u,v)$ for which $u\ge v$ and $v\ge p(u)$. If $u\le b$ or $u>a$, then $p(u)=u$, so $u=v$. If $b<u\le a$, then $p(u)=b$, so $b\le v\le u$. This gives precisely the family \eqref{eq:Omegaab}. It also proves (i) and (ii) without a separate closure calculation.

The initial identity-agreement interval for $p$ is $[0,b]$, while $q$ agrees with the identity everywhere. Theorem~\ref{thm:firstdeparture} proves (iii). Finally, the displayed description of the topologies gives proper inclusion for $a_1<a_2$. If $t\in(a_1,a_2]$, then $h_{xy}^{a_1,b}(t)=t$ but $h_{xy}^{a_2,b}(t)=b$. Choosing $\beta\in(b,t]$ and applying \eqref{eq:ordermatrix} separates the two membership-point preorders. Thus they also have distinct kernels.
\par The example shows explicitly that the graded relation cannot recover the membership-point preorder, and that identity agreement may resume above the first departure.
\end{example}

\begin{corollary} \label{thm:gradedfibrecard}\label{cor:gradedrigidity}
Fix $X=\{x,y\}$ and $(b,c)\in[0,1]^2$. Among stratified Alexandrov Chang fuzzy topologies with $r(x,y)=b$ and $r(y,x)=c$, there is exactly one topology if $(b,c)=(1,1)$, namely the topology of all constants. Otherwise there are exactly continuum many topologies and exactly continuum many distinct membership-point specialization preorders.
\end{corollary}
\begin{proof}
Every reflexive relation on two points is minimum-transitive: compositions with a repeated adjacent support reproduce the corresponding off-diagonal entry, and compositions returning to the starting support are bounded by $1$. Apply Theorem~\ref{thm:generalfibre}(ii)--(iv), and use Corollary~\ref{cor:orderrec} to identify distinct Alexandrov topologies with distinct membership preorders.
\end{proof}

\begin{theorem} \label{thm:fibreorder}
Let $X=\{x,y\}$ and fix $b,c\in[0,1]$. For a join-preserving map $f:[0,1]\to[0,1]$ with $f\le\operatorname{id}$, put
\[
 \delta(f)=\inf\bigl(\{t>0:f(t)<t\}\cup\{1\}\bigr),\qquad
 \mathcal J_d=\{f\in\operatorname{Sup}([0,1]): f\le\operatorname{id},\ \delta(f)=d\}.
\]
Let $\mathfrak A_{b,c}$ denote the collection of stratified Alexandrov Chang fuzzy topologies with $r(x,y)=b$ and $r(y,x)=c$, ordered by inclusion. Then:
\begin{enumerate}[label=\textup{(\roman*)}]
\item The association $(p,q)\mapsto\Omega_{p,q}$ of Theorem~\ref{thm:twopointclass} is an order anti-isomorphism
\[
 (\mathcal J_b\times\mathcal J_c,\le_{\mathrm{pt}})
 \ \cong^{\mathrm{op}}\ (\mathfrak A_{b,c},\subseteq),
\]
where the product carries the pointwise order. In particular,
\[
 \Omega_{p,q}\subseteq\Omega_{p',q'}
 \quad\Longleftrightarrow\quad p\ge p'\ \text{and}\ q\ge q'.
\]
\item Each $\mathcal J_d$ has a least element $m_d(t)=\min\{t,d\}$; correspondingly, $\mathfrak A_{b,c}$ has the greatest element $\Omega_r$, with $r(x,y)=b$, $r(y,x)=c$.
\item If $b=c=1$, the fiber $\mathfrak A_{b,c}$ is a singleton. Otherwise it has no inclusion-minimal element (and hence no least element), by Theorem~\ref{thm:generalfibre}.
\end{enumerate}
\end{theorem}
\begin{proof}
By Theorem~\ref{thm:twopointclass} and the first-departure formula, the indicated map is a bijection. For any two join-preserving closure operators $K,K'$ on $[0,1]^X$, their fixed-point families satisfy
\[
 \Fix(K)\subseteq\Fix(K')\quad\Longleftrightarrow\quad K\ge K'
\]
(pointwise): if $\Fix(K)\subseteq\Fix(K')$, then $K(A)\in\Fix(K')$ and $A\le K(A)$ imply $K'(A)\le K(A)$; conversely if $K\ge K'$, a fixed point of $K$ must also be fixed by $K'$ by extensivity. In the two-point formula for $K_{p,q}$, the inequality $K_{p,q}\ge K_{p',q'}$ is equivalent to $p\ge p'$ and $q\ge q'$ by evaluating at fuzzy points. This proves (i).

By Lemma~\ref{lem:thresholdbound}, every $f\in\mathcal J_d$ satisfies $\min\{t,d\}\le f(t)$. The map $m_d$ belongs to $\mathcal J_d$, so it is least. Its associated topology is $\Omega_r$ as in Theorem~\ref{thm:gradedmax}, proving (ii).

The singleton and absence-of-minimal-elements claims in (iii) follow directly from Theorem~\ref{thm:generalfibre}. Equivalently, $\mathcal J_1=\{\operatorname{id}\}$, whereas every $f\in\mathcal J_d$ with $d<1$ has a strictly larger member $f\vee(\lambda\operatorname{id})$ in the same class for a suitable $\lambda<1$, by the argument of Lemma~\ref{lem:enlarge}.
\end{proof}

\begin{remark}[Scope of the two-point classification]
Lai--Zhang \cite[Example~4.24]{LaiZhang2006} already exhibit
nonuniqueness using a topology generated by constants and a
truncated characteristic function. On two points, this is the
$b=0$ case of \eqref{eq:Omegaab}. Neither that example nor the
general scalar upper-set and matrix constructions is claimed as new.
Theorem~\ref{thm:fibreorder} instead identifies the entire inclusion
order of every fixed two-point Alexandrov fiber as the opposite of
a product of first-departure classes. Theorem~\ref{thm:generalfibre}
gives the complementary uniform descending-chain statement below
every stratified realization.
\end{remark}

\section{Kernel completions and quasi-specialization}\label{sec:families}
\subsection{Open families with prescribed membership information}
Theorem~\ref{thm:matrixrep} determines a unique Alexandrov topology from an admissible threshold matrix. General Chang fuzzy topologies are not necessarily Alexandrov, so the same matrix can arise from different open families. The following characterization states exactly how this occurs.

\begin{theorem} \label{thm:fibre}
Let $C\subseteq[0,1]^X$ be a stratified Alexandrov Chang fuzzy topology. A stratified Chang fuzzy topology $\Omega$ satisfies $\Fix(K_\Omega)=C$ if and only if
\begin{enumerate}[label=\textup{(\roman*)}]
\item $\Omega\subseteq C$;
\item every $D\in C$ can be written as $D=\bigwedge\mathcal U$ for some family $\mathcal U\subseteq\Omega$.
\end{enumerate}
Therefore the stratified fuzzy topologies with this membership-point preorder are exactly the join-closed, finite-meet-closed, constant-containing subfamilies of $C$ that are meet-dense in $C$.
\end{theorem}
\begin{proof}
If $\Fix(K_\Omega)=C$, then every open $U\in\Omega$ is fixed, so (i) holds. Every $D\in C$ satisfies $D=K_\Omega(D)$, the meet of all opens above it, giving (ii). Conversely, suppose (i) and (ii). Every fixed point of $K_\Omega$ is an arbitrary meet of members of $\Omega$, and hence belongs to $C$ by its Alexandrov property. By (ii), for $D=\bigwedge\mathcal U$ we have $D\le U$ for every $U\in\mathcal U$, so $D\le K_\Omega(D)\le\bigwedge\mathcal U=D$. Hence $\Fix(K_\Omega)=C$.
\end{proof}

This characterization also clarifies why the quasi-point specialization order carries different information from the membership order: the latter sees only the completion, while quasi-coincidence tests strict inequalities against the original open family.

\subsection{Recovery of membership information from quasi-specialization}
\begin{theorem} \label{thm:quasireconstruction}
For every Chang fuzzy topology on $X$, its quasi-point specialization preorder determines its membership-point specialization preorder and its open-kernel operator. No stratification assumption is needed.

More precisely, for $x,y\in X$ and $0\le s<1$ put
\[
 g_{xy}(s)=\inf\{U(y):U\in\Omega,\ U(x)>s\}.
\]
Then the preorder determines these numbers by
\begin{equation}\label{eq:quasig}
 g_{xy}(s)=\sup\bigl(\{b\in[0,1):x_{1-s}\preQ y_{1-b}\}\cup\{0\}\bigr),
\end{equation}
and the membership threshold profiles are recovered as
\begin{equation}\label{eq:quasitomembership}
 h_{xy}(0)=0,\qquad h_{xy}(t)=\sup_{0\le s<t}g_{xy}(s)\quad(0<t\le1).
\end{equation}
\end{theorem}
\begin{proof}
The set defining $g_{xy}(s)$ is nonempty because $1\in\Omega$. By the definition of quasi-specialization,
\[
 x_{1-s}\preQ y_{1-b}
 \quad\Longleftrightarrow\quad
 (U(x)>s\ \Longrightarrow\ U(y)>b)\text{ for all }U\in\Omega.
\]
Writing $g=g_{xy}(s)$, every $b\in[0,1)$ with $b<g$ satisfies this implication, and no $b>g$ does. Indeed, if $b>g$, some open in the defining infimum has $U(y)<b$. At $b=g$ the implication may or may not hold, depending on attainment of the infimum. Including $0$ handles the case $g=0$; if $g=1$ the supremum of $[0,1)$ is $1$. This proves \eqref{eq:quasig}.

For $0\le s<t\le1$ the inclusion of the open families is
\[
 \{U:U(x)\ge t\}\subseteq\{U:U(x)>s\}
 \subseteq\{U:U(x)\ge s\},
\]
where all opens are taken from $\Omega$. Consequently,
\[
 h_{xy}(s)\le g_{xy}(s)\le h_{xy}(t).
\]
The kernel preserves arbitrary joins by \eqref{eq:join}, even without stratification. Since $t=\sup_{0\le s<t}s$, it follows that
\[
 h_{xy}(t)=\sup_{0\le s<t}h_{xy}(s)
 \le\sup_{0\le s<t}g_{xy}(s)\le h_{xy}(t).
\]
This is \eqref{eq:quasitomembership}. The membership preorder is then recovered by $x_\alpha\preM y_\beta$ if and only if $\beta\le h_{xy}(\alpha)$, and Theorem~\ref{thm:kernelorder} gives the full kernel.
\end{proof}

\begin{corollary}\label{cor:quasialex}
Two Chang fuzzy topologies on the same labeled set with equal quasi-point specialization preorders have equal membership preorders and equal kernel completions. If both are Alexandrov, their original fuzzy topologies coincide.
\end{corollary}
\begin{proof}
Apply Theorem~\ref{thm:quasireconstruction} and use $\Omega=\Fix(K_\Omega)$ for Alexandrov topologies.
\end{proof}

\subsection{Finite-set reconstruction}

The Alexandrov assumption is not necessary for recovering the original open family from quasi-specialization when the support set is finite. The following argument uses finite meets in the Chang topology, rather than arbitrary meets in its kernel completion.

\begin{theorem} \label{thm:finitequasi}
Let $X$ be finite and nonempty, and let $\Omega$ be any Chang fuzzy topology on $X$. Then
\begin{equation}\label{eq:finitequasi}
 \Omega=\{A\in[0,1]^X:\Qq(A)\text{ is upward closed for }\preQ\}.
\end{equation}
In particular, two Chang fuzzy topologies on the same finite labeled set have equal quasi-specialization preorders if and only if they coincide. No stratification or Alexandrov assumption is needed.
\end{theorem}
\begin{proof}
For $A\in\Omega$, upward closure of $\Qq(A)$ follows from the definition of the quasi-specialization preorder. Conversely, suppose that $\Qq(A)$ is upward closed. Fix $x\in X$ and $0\le s<A(x)$. For each $y\in X$ there is an open $U_y\in\Omega$ with
\[
 U_y(x)>s,\qquad U_y(y)\le A(y).
\]
If $A(y)=1$, take $U_y=1$. Otherwise, if no such open exists, then every open satisfies $U(x)>s\Rightarrow U(y)>A(y)$. This means $x_{1-s}\preQ y_{1-A(y)}$. But $x_{1-s}\in\Qq(A)$ and $y_{1-A(y)}\notin\Qq(A)$, a contradiction.

Because $X$ is finite, $V_{x,s}=\bigwedge_{y\in X}U_y$ is open, satisfies $V_{x,s}\le A$, and has $V_{x,s}(x)>s$. Taking the join of these opens over all $x$ and $0\le s<A(x)$ gives $A$: each summand is below $A$, and the supremum at every coordinate with $A(x)>0$ is at least $\sup_{s<A(x)}s=A(x)$. The zero case is the empty join. Thus $A\in\Omega$.
\end{proof}

The formula reconstructs the fuzzy open family, not the entire family of upward closed subsets of the point preorder: only sets of the form $\Qq(A)$ enter \eqref{eq:finitequasi}. In particular, it does not assert that the quasi-point topology itself is Alexandrov.

\subsection{Greatest quasi-specialization realizations and finite restrictions}
For a Chang fuzzy topology $\Omega$ on $X$, write $R_q=\preQ^{\Omega}$ and define
\begin{equation}\label{eq:qmax}
 \Omega_q^*=\{A\in[0,1]^X:\Qq(A)\in\Up(R_q)\}.
\end{equation}
Here $\Up(R_q)$ denotes the upward closed subsets of the labeled fuzzy-point universe. For a nonempty $F\subseteq X$, put $\Omega|_F=\{U|_F:U\in\Omega\}$. Restriction preserves arbitrary joins, finite meets, and the zero and unit functions, so $\Omega|_F$ is a Chang topology.

\begin{theorem} \label{thm:qmax}
Let $\Omega$ be a Chang fuzzy topology on a nonempty set $X$. Then:
\begin{enumerate}[label=\textup{(\roman*)}]
\item $\Omega_q^*$ is a Chang fuzzy topology containing $\Omega$ and has quasi-specialization preorder $R_q$. It is the greatest Chang topology with this preorder. If $\Omega$ is stratified, so is $\Omega_q^*$.
\item The following descriptions coincide:
\begin{align}
 \Omega_q^*&=\{A\in[0,1]^X:A|_F\in\Omega|_F
      \text{ for every finite nonempty }F\subseteq X\}\label{eq:qfinite}\\
 &=\{A\in[0,1]^X:A|_F\in\Omega|_F
      \text{ for every }F\subseteq X\text{ with }1\le |F|\le2\}.
 \label{eq:qpair}
\end{align}
Equivalently, $\Omega_q^*$ is the closure of $\Omega$ in the product space $[0,1]_{\mathrm{disc}}^X$, with the discrete topology on each copy of the grade interval.
\item $K_{\Omega_q^*}=K_\Omega$, and
\begin{equation}\label{eq:greatesthierarchy}
 \Omega\subseteq\Omega_q^*\subseteq\Fix(K_\Omega)
 \subseteq\Omega_{r_\Omega}.
\end{equation}
The last three topologies are the greatest realizations of the quasi-point preorder, the membership-point preorder, and the G\"odel relation of $\Omega$, respectively, among all Chang topologies on $X$. If $X$ is finite or $\Omega$ is Alexandrov, then $\Omega_q^*=\Omega$.
\end{enumerate}
\end{theorem}
\begin{proof}
(i) The empty set and $E$ are upward closed for $R_q$. By \eqref{eq:Q-ops}, arbitrary joins and finite meets of members of $\Omega_q^*$ remain in that family. Every $U\in\Omega$ has an upward closed quasi-spectrum, so $\Omega\subseteq\Omega_q^*$. This inclusion implies $\preQ^{\Omega_q^*}\subseteq R_q$. Conversely, if $pR_q q$, every $A\in\Omega_q^*$ satisfies $p\in\Qq(A)\Rightarrow q\in\Qq(A)$, so $R_q\subseteq\preQ^{\Omega_q^*}$. Thus the preorders agree. If another Chang topology $\Lambda$ has preorder $R_q$, every $V\in\Lambda$ has an upward closed quasi-spectrum for $R_q$, whence $\Lambda\subseteq\Omega_q^*$. Stratification follows from inclusion of $\Omega$.

(ii) For any nonempty $F\subseteq X$, the quasi-specialization preorder of $\Omega|_F$ is the restriction of $R_q$ to the fuzzy points supported in $F$: both preorders test precisely the same inequalities on pairs of coordinates in $F$. Hence, for $A\in\Omega_q^*$, the quasi-spectrum of $A|_F$ is upward closed for the restricted preorder. When $F$ is finite, Theorem~\ref{thm:finitequasi} gives $A|_F\in\Omega|_F$. This proves the inclusion from $\Omega_q^*$ into the right side of \eqref{eq:qfinite}; that family is in turn contained in the right side of \eqref{eq:qpair}.

Conversely, assume $A$ belongs to the right side of \eqref{eq:qpair}. If $x_\alpha R_q y_\beta$ and $x_\alpha\in\Qq(A)$, take $F=\{x,y\}$ (a singleton if $x=y$). There is $U\in\Omega$ with $U|_F=A|_F$. Then
\[
 A(x)>1-\alpha\ \Longrightarrow\ U(x)>1-\alpha
 \ \Longrightarrow\ U(y)>1-\beta
 \ \Longrightarrow\ A(y)>1-\beta.
\]
Thus $\Qq(A)$ is upward closed and $A\in\Omega_q^*$. This proves both equalities. In the product of discrete grade intervals, a basic neighborhood of $A$ specifies its exact values on finitely many coordinates. It meets $\Omega$ precisely when those values occur as the restriction of some member of $\Omega$. This is exactly \eqref{eq:qfinite}, proving the closure description.

(iii) By (i) and Theorem~\ref{thm:quasireconstruction}, $K_{\Omega_q^*}=K_\Omega$. Every open is fixed by its own kernel, proving the middle inclusion in \eqref{eq:greatesthierarchy}. The completion $\Fix(K_\Omega)$ itself has kernel $K_\Omega$: for each $A$, $K_\Omega(A)$ is the least fixed point above $A$. Any topology with that kernel consists of fixed points of $K_\Omega$, so the completion is its greatest membership-preorder realization by Theorem~\ref{thm:kernelorder}. By Theorem~\ref{thm:gradedprofile}, it also has G\"odel relation $r_\Omega$. That relation is reflexive and minimum-transitive (see Section~\ref{sec:loss}), so Theorem~\ref{thm:gradedmax} gives the final inclusion and greatest graded realization. Part (i) supplies the greatest quasi-realization assertion. If $X$ is finite, \eqref{eq:qfinite} with $F=X$ yields $\Omega_q^*=\Omega$. If $\Omega$ is Alexandrov, $\Omega=\Fix(K_\Omega)$, and \eqref{eq:greatesthierarchy} yields the same equality.
\end{proof}

\begin{corollary} \label{cor:qfinite}
For two Chang fuzzy topologies $\Omega$ and $\Lambda$ on the same nonempty labeled set $X$, the following are equivalent:
\begin{enumerate}[label=\textup{(\roman*)}]
\item Their quasi-specialization preorders coincide.
\item $\Omega_q^*=\Lambda_q^*$.
\item $\Omega|_F=\Lambda|_F$ for every finite nonempty $F\subseteq X$.
\item $\Omega|_F=\Lambda|_F$ for every $F\subseteq X$ with $1\le |F|\le2$.
\end{enumerate}
In particular, $\Omega$ is the greatest topology with its quasi-specialization preorder if and only if it is closed in $[0,1]_{\mathrm{disc}}^X$.
\end{corollary}
\begin{proof}
The implication (i)$\Rightarrow$(ii) follows from \eqref{eq:qmax}, and the converse follows from Theorem~\ref{thm:qmax}(i). For every finite nonempty $F$, \eqref{eq:qfinite} and $\Omega\subseteq\Omega_q^*$ give $\Omega_q^*|_F=\Omega|_F$, and likewise for $\Lambda$. This proves (ii)$\Rightarrow$(iii). Clearly (iii)$\Rightarrow$(iv), and \eqref{eq:qpair} gives (iv)$\Rightarrow$(ii). The final claim follows from the closure description and greatest-realization assertion of Theorem~\ref{thm:qmax}.
\end{proof}

The grade coordinates in the product topology are discrete:
restrictions must match exactly, not approximately. The
finite-product two-coordinate projection principle is classical
(Topkis \cite[Theorem~1]{Topkis1976}; compare Baker--Pixley).
Theorem~\ref{thm:qmax} identifies the greatest quasi-specialization
realization with this finite-coordinate closure, even for infinite
$X$. Its quasi-spectra have the form $\Qq(A)$; the theorem does not
identify all upper sets of the fuzzy-point preorder with fuzzy
open sets.

\subsection{Cardinality on finite support sets}
\begin{lemma} \label{lem:countquasi}
On a fixed finite nonempty set $X$, there are at most $2^{\aleph_0}$ distinct quasi-specialization preorders arising from Chang fuzzy topologies.
\end{lemma}
\begin{proof}
Fix an ordered pair $x,y\in X$ and set
\[
 g(s)=\inf\{U(y):U\in\Omega,\ U(x)>s\},\qquad 0\le s<1.
\]
The set of opens in this infimum is nonempty, since $1\in\Omega$. The function $g$ is nondecreasing. More precisely, for every $s<1$,
\begin{equation}\label{eq:grightcount}
 g(s)=\inf_{\substack{u\in\mathbb Q\\s<u<1}}g(u).
\end{equation}
Indeed, every open with $U(x)>s$ has $U(x)>u$ for some rational $u\in(s,U(x))$; hence the family defining $g(s)$ is the union of the families defining $g(u)$ over those rational $u$. Taking infima proves \eqref{eq:grightcount}. Thus the values of $g$ on the countable set $\mathbb Q\cap[0,1)$ determine all its values, yielding at most continuum many functions $g$.

The quasi-specialization implication for $x_{1-s}$ and $y_{1-b}$, $0\le b<1$, is
\[
 x_{1-s}\preQ y_{1-b}
 \quad\Longleftrightarrow\quad
 U(y)>b\text{ for every }U\in\Omega\text{ with }U(x)>s.
\]
For fixed $g$, the implication holds if $b<g(s)$ and fails if
$b>g(s)$. When $b=g(s)<1$, the strict inequality in the definition
of quasi-coincidence matters: the implication holds exactly when
the defining infimum is not attained. Thus $g$ alone does not
determine the quasi-point preorder. For example, in
Example~\ref{ex:samekernel}, with $x=y$ and $0\le s<c$, the two
families $\{0\}\cup(c,1]$ and $\{0\}\cup[c,1]$ have the same value
$g(s)=c$, but the second family attains the infimum and the first
does not. Record this additional boundary information by
\[
 D=\{s\in[0,1):\exists U\in\Omega,\ U(x)>s,\ U(y)=g(s)\}.
\]

We show that, for each fixed $g$, the set $D$ has at most continuum many possible forms.

For a value $v$ in the range of $g$, its level set $L_v=\{s\in[0,1):g(s)=v\}$ is an interval (possibly a singleton), by monotonicity. If $s\in D$ and $U$ witnesses the defining property, then for each $u$ with $s\le u<U(x)$, monotonicity and this same witness give
\[
 g(s)\le g(u)\le U(y)=g(s),
\]
so $[s,U(x))\subseteq L_{g(s)}$. In particular, $L_{g(s)}$ contains a nondegenerate interval. There are at most countably many distinct nondegenerate level sets: their interiors are disjoint nonempty intervals, and each contains a distinct rational number.

For each such level set $L_v$, the subset $D\cap L_v$ is an initial segment of $L_v$. Indeed, if $s_1<s_2$ lie in $L_v$ and $s_2\in D$ with witness $U$, then $U(x)>s_2>s_1$ and $U(y)=g(s_2)=g(s_1)=v$, so $s_1\in D$. An initial segment of a real interval is specified by a cut point and, where applicable, whether that endpoint is included; there are at most $2^{\aleph_0}$ such segments. The set $D$ is empty on all singleton level sets, by the previous paragraph. Since there are only countably many nondegenerate level sets, the number of possible $D$ for fixed $g$ is at most
\[
 (2^{\aleph_0})^{\aleph_0}=2^{\aleph_0}.
\]
Hence the pair $(g,D)$ determines the quasi-specialization implications for this support pair and admits at most continuum many possibilities. When $g(s)=0$, the endpoint $b=0$ is treated by $D$; when $g(s)=1$, no boundary occurs because $b<1$. Finally, $X$ has finitely many ordered support pairs, so the finite product of their cardinal bounds is again $2^{\aleph_0}$.
\end{proof}

\begin{theorem} \label{thm:finiteallcard}
Let $X$ be finite and nonempty and let $r:X\times X\to[0,1]$ be reflexive and minimum-transitive. Then the family $\mathfrak T_r$ of all stratified Chang realizations satisfies
\begin{equation}\label{eq:finiteallcard}
 |\mathfrak T_r|=\begin{cases}
 1,&r\equiv1,\\
 2^{\aleph_0},&r\not\equiv1.
 \end{cases}
\end{equation}
In fact, on every fixed finite nonempty set there are exactly $2^{\aleph_0}$ Chang fuzzy topologies when stratification is not required.
\end{theorem}
\begin{proof}
Theorem~\ref{thm:finitequasi} identifies each Chang topology on $X$ uniquely by its quasi-specialization preorder. Lemma~\ref{lem:countquasi} therefore bounds the number of all Chang topologies by $2^{\aleph_0}$. The singleton case and the continuum lower bound for nontrivial $r$ follow from Theorem~\ref{thm:generalfibre}. Combining the bounds proves \eqref{eq:finiteallcard}. For the unrestricted count, each $t\in(0,1)$ gives the three-element Chang topology $\{0,\underline t,1\}$. These families are distinct because $X$ is nonempty, giving the matching continuum lower bound for all Chang topologies. On a singleton, this unrestricted count should not be confused with the unique stratified topology.
\end{proof}

The upper bound uses both the monotone threshold data and their boundary-attainment information. Counting only the membership threshold matrices would not suffice here, since distinct non-Alexandrov open families can have the same kernel. The next example makes this distinction explicit even on two points.

\subsection{Limitations and the reconstruction hierarchy}
\begin{example} \label{ex:quasifibre}
Fix $0<b<a<1$ and let $X=\{x,y\}$; write a fuzzy subset as its coordinate pair $(u,v)$. Consider
\begin{align*}
 C&=\{\underline t:t\in[0,1]\}
 \cup\{(u,v):b\le v<u\le a\},\\
 \Omega^\circ&=\{\underline t:t\in[0,1]\}
 \cup\{(u,v):b<v<u\le a\}.
\end{align*}
The family $C$ is the stratified Alexandrov topology $\Omega_{a,b}$ of Example~\ref{thm:continuum}. The family $\Omega^\circ$ is a stratified Chang fuzzy topology: it contains every constant and is closed under finite coordinatewise minima and arbitrary coordinatewise suprema. To verify closure, note that a nonconstant result of either operation must still have $v>b$ and $u\le a$; the alternative is a constant pair. (For an arbitrary join, if any constant exceeds $a$, the resulting join is constant.)

Clearly $\Omega^\circ\subsetneq C$. Every missing pair $(u,b)$ with $b<u\le a$ is an arbitrary meet of elements of $\Omega^\circ$, for example
\[
 (u,b)=\bigwedge_{n\ge1}
 \left(u,b+\frac{u-b}{n+1}\right).
\]
Consequently $\Omega^\circ$ is meet-dense in $C$; since $C$ is Alexandrov, Theorem~\ref{thm:fibre} gives
\[
 \Fix(K_{\Omega^\circ})=C=\Fix(K_C),
 \qquad K_{\Omega^\circ}=K_C.
\]
Theorem~\ref{thm:kernelorder} shows that the \emph{membership} specialization preorders coincide, and the G\"odel relations coincide as well.

Their \emph{quasi-point} orders differ. Take $b=1/4$, $a=3/4$, $p=x_{1/2}$, and $q=y_{3/4}$. Every $U\in\Omega^\circ$ satisfying $U(x)>1/2$ also satisfies $U(y)>1/4$: this is immediate for constants and follows from $v>b$ for nonconstant opens. Thus $p\preccurlyeq_q^{\Omega^\circ}q$. But $U=(0.6,1/4)\in C$ satisfies $U(x)>1/2$ and $U(y)\not>1/4$, so $p\not\preccurlyeq_q^Cq$. In particular, the membership specialization preorder does not determine quasi-specialization \emph{even among stratified Chang topologies}. By finite-set reconstruction, $(\Omega^\circ)_q^*=\Omega^\circ\subsetneq C=\Fix(K_{\Omega^\circ})$, so the middle inclusion of \eqref{eq:greatesthierarchy} can be strict.
\end{example}

\begin{example} \label{ex:quasinotop}
Let $X$ be any infinite set, let $\tau_1$ be the cofinite topology and $\tau_2$ the discrete topology, and define
\[
 \Omega_i=\{U:X\to[0,1]: U\text{ is lower semicontinuous for }\tau_i\},\qquad i=1,2.
\]
Both families are stratified Chang fuzzy topologies: the strict cuts $\{x:U(x)>s\}$ of an arbitrary join are unions of the corresponding strict cuts, and those of a finite meet are finite intersections. Constants are lower semicontinuous. Clearly $\Omega_1\subseteq\Omega_2$; the inclusion is strict because the characteristic function of a singleton belongs to $\Omega_2$ but not to $\Omega_1$.

Nevertheless, both quasi-point specialization preorders are given by
\begin{equation}\label{eq:quasidisc}
 x_\alpha\preQ y_\beta\quad\Longleftrightarrow\quad x=y\text{ and }\alpha\le\beta.
\end{equation}
Indeed, if $x\ne y$, the characteristic function of $X\setminus\{y\}$ belongs to both topologies, takes value $1$ at $x$ and $0$ at $y$, and is quasi-coincident with $x_\alpha$ but not $y_\beta$. For $x=y$ and $\alpha\le\beta$, the implication $U(x)>1-\alpha\Rightarrow U(x)>1-\beta$ holds for every open. If $\alpha>\beta$, the constant with value $1-\beta$ is quasi-coincident with $x_\alpha$ but not with $x_\beta$. This proves \eqref{eq:quasidisc}. Thus quasi-specialization need not recover the original fuzzy open family, even with stratification. In this example $(\Omega_1)_q^*=\Omega_2=[0,1]^X$: every fuzzy set has an upward closed quasi-spectrum for \eqref{eq:quasidisc}. Equivalently, every finite restriction of any fuzzy set extends to a lower semicontinuous function for the cofinite topology.
\end{example}

\begin{proposition} \label{prop:contrast}
For arbitrary Chang fuzzy topologies on a fixed labeled fuzzy-point universe:
\begin{enumerate}[label=\textup{(\roman*)}]
\item either full classical point topology $\tau_m$ or $\tau_q$ determines the original fuzzy topology;
\item the quasi-point preorder determines the original fuzzy topology on finite support sets, but need not do so on any infinite support set, even under stratification;
\item the quasi-point preorder determines the membership-point preorder, which is equivalent to the open-kernel operator and determines exactly the Alexandrov completion. The membership preorder need not determine the quasi-point preorder, even under stratification;
\item the membership preorder determines the G\"odel relation without stratification, but the converse fails even for stratified Alexandrov topologies.
\end{enumerate}
\end{proposition}
\begin{proof}
Part (i) follows from Theorem~\ref{thm:mrecovery} and Proposition~\ref{prop:Qframe}. Part (ii) follows from Theorem~\ref{thm:finitequasi} and Example~\ref{ex:quasinotop}. For (iii), apply Theorems~\ref{thm:kernelorder} and~\ref{thm:quasireconstruction} and Example~\ref{ex:quasifibre}. For (iv), apply Theorem~\ref{thm:gradedprofile} and the distinct Alexandrov realizations in Theorem~\ref{thm:generalfibre}, or the explicit two-point family of Example~\ref{thm:continuum}.
\end{proof}

For arbitrary Chang fuzzy topologies, the resulting reconstruction hierarchy is
\[
 \Omega\ \longrightarrow\ \preQ\ \longrightarrow\ \preM
 \ \longleftrightarrow\ K_\Omega\ \longrightarrow\ r_\Omega.
\]
Every one-way arrow is strict in general: its reverse determination fails, even among stratified topologies. The examples witnessing the failures, in order from left to right, are Examples~\ref{ex:quasinotop} and~\ref{ex:quasifibre} and Example~\ref{thm:continuum}. In the Alexandrov class the first two arrows are reversible, since either point preorder recovers the kernel and the kernel recovers the original topology. On finite support sets the first arrow is reversible even outside the Alexandrov class, but the second can still be strict by Example~\ref{ex:quasifibre}. The last arrow is not reversible in the Alexandrov class, even under stratification. Either full point topology determines $\Omega$ independently of the Alexandrov assumption. The hierarchy of greatest realizations \eqref{eq:greatesthierarchy} runs in the opposite direction of information loss: each of its three inclusions can be strict, as witnessed respectively by Examples~\ref{ex:quasinotop}, \ref{ex:quasifibre}, and the non-scalar Alexandrov realizations of Section~\ref{sec:twopoint}.

\section{Maps, point-order preservation, and fuzzy continuity}\label{sec:order-maps}
A map $f:X\to Y$ has a grade-preserving lift $\widehat f(x_\alpha)=f(x)_\alpha$. Standard fuzzy continuity implies continuity of the lift for both point topologies, as recalled later. Requiring only preservation of the membership specialization order is weaker; it has an exact description using the kernel.

\begin{theorem} \label{thm:laxiff}
Let $(X,\Omega_X)$ and $(Y,\Omega_Y)$ be Chang fuzzy topological spaces, and let $f:X\to Y$ be any map. The following are equivalent:
\begin{enumerate}[label=\textup{(\roman*)}]
\item $\widehat f$ preserves the membership-point specialization preorder;
\item for every fuzzy subset $A$ of $Y$,
\begin{equation}\label{eq:lax}
 K_X(f^*A)\le f^*K_Y(A),\qquad f^*A=A\circ f.
\end{equation}
\end{enumerate}
These are also equivalent to the explicit inequalities
\[
 h^X_{xz}(t)\le h^Y_{f(x)f(z)}(t)
 \qquad(x,z\in X,\ t\in[0,1]).
\]
\end{theorem}
\begin{proof}
Suppose (i). For $x,z\in X$ and $t>0$, order preservation and \eqref{eq:ordermatrix}, whose equivalence with kernel membership does not require stratification, give
\[
 K_X(x_t)(z)\le K_Y(f(x)_t)(f(z)).
\]
Writing $f^*A=\bigvee_{x\in X}x_{A(f(x))}$ and applying arbitrary-join preservation yields
\[
 K_X(f^*A)(z)
 =\sup_x K_X(x_{A(f(x))})(z)
 \le\sup_xK_Y(f(x)_{A(f(x))})(f(z))
 \le K_Y(A)(f(z)).
\]
Thus (ii) holds. Conversely, suppose (ii). Apply it to the fuzzy point $A=f(x)_t$, and note that $x_t\le f^*A$. Monotonicity gives
\[
 K_X(x_t)(z)\le K_X(f^*A)(z)
 \le K_Y(f(x)_t)(f(z)),
\]
which is exactly the order-preservation condition. The last assertion is the same inequality in the notation of Section~\ref{sec:profiles}.
\end{proof}

\begin{corollary} \label{cor:alexcont}
If $\Omega_X$ is Alexandrov, then $f:(X,\Omega_X)\to(Y,\Omega_Y)$ is fuzzy-continuous if and only if its grade-preserving lift preserves membership-point specialization. If both topologies are Alexandrov, the point-order comparison completely characterizes fuzzy-continuous maps between them.
\end{corollary}
\begin{proof}
Continuity implies \eqref{eq:lax}: inverse images of open supersets of $A$ provide open supersets of $f^*A$, so $K_X(f^*A)\le f^*K_Y(A)$. Conversely, suppose \eqref{eq:lax}. For $V\in\Omega_Y$ we have $K_Y(V)=V$, so
\[
 f^*V\le K_X(f^*V)\le f^*K_Y(V)=f^*V.
\]
Thus $f^*V$ is fixed by $K_X$. Since $\Omega_X$ is Alexandrov, $\Omega_X=\Fix(K_X)$, so $f^*V$ is open. This proves continuity.
\end{proof}

The next example shows why the Alexandrov assumption matters.

\begin{example}\label{ex:nonalexmap}
Take a non-Alexandrov Chang fuzzy topology $\Omega$ on a set $X$, for instance $X=\{x\}$ with $\Omega=\{0\}\cup(c,1]$. Let $\Lambda=\Fix(K_\Omega)$ and $f=\operatorname{id}_X:(X,\Omega)\to(X,\Lambda)$. Both topologies have the same kernel, so $\widehat f$ preserves their membership-point orders. But $f$ is not fuzzy-continuous, since the constant $c$ is open in $\Lambda$ and is not open in $\Omega$.
\end{example}

\begin{corollary} \label{cor:completionmaps}
For every set map $f:X\to Y$, its grade-preserving lift preserves membership specialization if and only if
\[
 f:(X,\Fix(K_X))\longrightarrow(Y,\Fix(K_Y))
\]
is fuzzy-continuous.
\end{corollary}
\begin{proof}
The two completions have the same kernels and hence the same membership specialization orders as the original spaces. Both completions are Alexandrov. Apply Corollary~\ref{cor:alexcont}.
\end{proof}

\section{Conclusion}
For every Chang fuzzy topology the quasi-point specialization preorder determines the membership preorder and its kernel, and the latter determines G\"odel specialization. Each one-way reduction loses information in general, even for stratified spaces. Either full point preorder recovers an Alexandrov fuzzy topology; quasi-specialization also recovers every Chang fuzzy topology on a finite support set. More generally, its greatest realizing topology consists exactly of the fuzzy sets that match original opens on each one- or two-point subset, equivalently on each finite subset. It is the closure of the original family in the product of discrete grade intervals. Together with the membership completion and the scalar upper-set topology, this yields a nested hierarchy of greatest realizations for the three specialization data. The graded relation generally does not recover the full profiles: under stratification it retains only the endpoint of each profile's initial interval of identity agreement.

For every prescribed reflexive, minimum-transitive fuzzy preorder, the stratified Chang realization fiber has a greatest topology. It is a singleton exactly when the relation is identically $1$; otherwise every realization strictly contains a continuum-sized descending chain with no lower bound in the fiber. The restriction construction preserves the Alexandrov property. Nontrivial fibers have exactly continuum cardinality for all stratified Chang topologies on finite sets and for Alexandrov topologies on at most countable sets. The total number of Chang topologies without a stratification requirement on every fixed finite nonempty set is also exactly continuum. On two points the product representation describes the entire Alexandrov fiber order. Scalar compatibility is exactly the extra kernel condition that reduces the profiles to $\min(t,r(x,y))$. These results distinguish full threshold matrices and strict-boundary information from the established scalar fuzzy-preorder representation. Further questions include the exact cardinality of general stratified fibers on infinite sets, the size and order structure of Alexandrov fibers on uncountable sets, and reconstruction from additional graded observations.

\appendix
\section{Supplementary properties of the point spaces}\label{sec:other}
The remaining results record elementary consequences of the two point-space constructions. They are included for completeness and are not proposed as the principal contributions.
\subsection{Vertical subspaces and the projection to the underlying set}\label{sec:projection}
The topology of the individual membership-level fibers is simple in a stratified fuzzy space. It also helps explain the difference between the two point topologies. Let $E_x=\{x_\alpha:0<\alpha\le1\}$.

\begin{proposition} \label{prop:verticaltop}
Suppose $\Omega$ is stratified. Under the identification $E_x\cong(0,1]$, the subspace topology of $\tauM$ consists of all downward closed subsets of $(0,1]$, and the subspace topology of $\tauQ$ consists of the empty set and all intervals $(a,1]$, $0\le a<1$. In particular,
\[
 x_\alpha\preM x_\beta\iff\alpha\ge\beta,
 \qquad
 x_\alpha\preQ x_\beta\iff\alpha\le\beta.
\]
The restriction of the topology $\tauM\vee\tauQ$ to $E_x$ is generated by half-open intervals $(a,b]$ with $0\le a<b\le1$.
\end{proposition}
\begin{proof}
All vertical sections of membership opens are downward closed, by Lemma~\ref{lem:fibers}. Because all constant fuzzy sets are open, the sets $\Pp(\underline c)\cap E_x=(0,c]$, for $0<c\le1$, are available. Arbitrary unions of these intervals give all downward closed subsets, including intervals $(0,c)$ obtained by unions over parameters below $c$.

Similarly, all quasi-open vertical sections are of the form $(a,1]$ or empty, and the constant fuzzy sets give every such section as $\Qq(\underline{1-a})\cap E_x$. The specialization formulas follow by testing which of these lower or upper sets contain a given membership level.

Finally, intersections of the basic lower sets $(0,b]$ and upper sets $(a,1]$ yield $(a,b]$. Every restriction of an intersection of a membership-open and a quasi-open set is a union of such intervals (and possibly intervals with an open upper endpoint, which are unions of half-open intervals). Thus these intervals form a base for the join of the restricted topologies.
\end{proof}

There is also a natural map forgetting the grade:
\[
 \pi:E\longrightarrow X,\qquad \pi(x_\alpha)=x.
\]
Put
\[
 \mathcal T_{\mathrm{cr}}=\{B\subseteq X:\chi_B\in\Omega\}.
\]
The collection $\mathcal T_{\mathrm{cr}}$ is a classical topology on $X$, because joins and finite meets of characteristic functions are the characteristic functions of unions and intersections.

\begin{theorem} \label{thm:projection}
For either point topology $\tauM$ or $\tauQ$, the map
\[
 \pi:(E,\tauM)\longrightarrow (X,\mathcal T_{\mathrm{cr}})
 \quad\text{or}\quad
 \pi:(E,\tauQ)\longrightarrow (X,\mathcal T_{\mathrm{cr}})
\]
is a continuous surjective quotient map. Moreover, the following conditions are equivalent:
\begin{enumerate}[label=\textup{(\roman*)}]
\item $\pi:(E,\tauM)\to(X,\mathcal T_{\mathrm{cr}})$ is open;
\item $\pi:(E,\tauQ)\to(X,\mathcal T_{\mathrm{cr}})$ is open;
\item $\chi_{\{x:U(x)>0\}}\in\Omega$ for every $U\in\Omega$.
\end{enumerate}
\end{theorem}
\begin{proof}
For any $B\subseteq X$,
\[
 \pi^{-1}(B)=B\times(0,1].
\]
If $B\in\mathcal T_{\mathrm{cr}}$, this set is simultaneously $\Pp(\chi_B)$ and $\Qq(\chi_B)$, hence open in both topologies. Conversely, if the preimage is open in either point topology, it belongs to that topology with full-or-empty vertical sections. By Theorem~\ref{thm:mrecovery} in the membership case, or Proposition~\ref{prop:Qframe} in the quasi case, it follows that $\chi_B\in\Omega$. Thus $\mathcal T_{\mathrm{cr}}$ is exactly the quotient topology under $\pi$, which proves the first assertion.

For each $U\in\Omega$, both $\pi(\Pp(U))$ and $\pi(\Qq(U))$ equal $\{x:U(x)>0\}$. If either projection is open, every such support belongs to $\mathcal T_{\mathrm{cr}}$, giving (iii). Conversely, if (iii) holds, the images of all basic membership opens are open. Since every membership-open set is a union of basic sets, $\pi$ is open on $(E,\tauM)$. The same argument for the quasi topology uses the sets $\Qq(U)$ as an open base, proving (ii).
\end{proof}

\begin{example} \label{ex:notopen}
Let $X=\{x,y\}$ and let $\Omega$ be the smallest stratified Chang fuzzy topology containing the fuzzy set $V=(0.6,0)$, where coordinates refer to $(x,y)$. Every open set $U\in\Omega$ satisfies
\[
 U(x)>U(y)\Longrightarrow U(x)\le0.6.
\]
Indeed, the displayed condition holds for constants and $V$, and it is preserved by arbitrary joins and finite meets, by elementary coordinatewise calculations. Therefore $\chi_{\{x\}}=(1,0)$ does not belong to $\Omega$, although $\{x\}$ is the support of $V$. The image of the open set $\Pp(V)$ (or $\Qq(V)$) under $\pi$ is $\{x\}$, which is not open in $\mathcal T_{\mathrm{cr}}$. Thus $\pi$ is a quotient map but not an open map.
\end{example}

\subsection{Continuity detected by either point topology}\label{sec:maps}
Let $f:X\to Y$ be a function between Chang fuzzy topological spaces $(X,\Omega_X)$ and $(Y,\Omega_Y)$. It is fuzzy-continuous if $V\circ f\in\Omega_X$ for every $V\in\Omega_Y$. Define its grade-preserving lift
\[
 \widehat f:E_X\to E_Y,\qquad \widehat f(x_\alpha)=f(x)_\alpha.
\]
The following gives the exact continuity tests for the two point topologies.

\begin{theorem} \label{thm:continuity}
For any map $f:X\to Y$, the following are equivalent:
\begin{enumerate}[label=\textup{(\roman*)}]
\item $f:(X,\Omega_X)\to(Y,\Omega_Y)$ is fuzzy-continuous;
\item $\widehat f:(E_X,\tau_m(\Omega_X))\to(E_Y,\tau_m(\Omega_Y))$ is continuous;
\item $\widehat f:(E_X,\tau_q(\Omega_X))\to(E_Y,\tau_q(\Omega_Y))$ is continuous.
\end{enumerate}
\end{theorem}
\begin{proof}
For every $V\in\Omega_Y$, the definitions give
\begin{equation}\label{eq:inverse}
 \widehat f^{-1}(\Pp(V))=\Pp(V\circ f),\qquad
 \widehat f^{-1}(\Qq(V))=\Qq(V\circ f).
\end{equation}
If $f$ is fuzzy-continuous, these inverse images are open in the corresponding point topologies. This proves (i)$\Rightarrow$(ii) and (i)$\Rightarrow$(iii).

Suppose (ii). For each $V\in\Omega_Y$, the first set in \eqref{eq:inverse} is open in $\tau_m(\Omega_X)$. Each of its vertical sections is either empty or $(0,V(f(x))]$, and hence is closed in $(0,1]$. Theorem~\ref{thm:mrecovery} implies $V\circ f\in\Omega_X$. Thus (ii)$\Rightarrow$(i).

Suppose (iii). By Proposition~\ref{prop:Qframe}, the second inverse image in \eqref{eq:inverse} is $\Qq(U)$ for some $U\in\Omega_X$. Injectivity of the map $\Qq$ gives $U=V\circ f$. Hence $f$ is fuzzy-continuous.
\end{proof}

\begin{remark}
Grade preservation is essential: an arbitrary continuous point-space map need not arise from a map of the underlying fuzzy spaces.
\end{remark}

\subsection{Fixed membership levels and the strict-cut distinction}\label{sec:cuts}
For $\alpha\in(0,1]$ let $E_\alpha=\{x_\alpha:x\in X\}$. For $U\in\Omega$ define the nonstrict and strict cuts
\[
 [U]_{\alpha}=\{x:U(x)\ge\alpha\},\qquad
 [U]^>_t=\{x:U(x)>t\} \quad(0\le t<1).
\]
The membership topology restricted to $E_\alpha$ corresponds to the topology generated by the nonstrict cuts $[U]_\alpha$. The quasi-point topology restricted to $E_\alpha$ corresponds to the strict-cut family $\{[U]^>_{1-\alpha}:U\in\Omega\}$, which is already a topology. These follow from the point conventions \cite{PuLiu1980}.

\begin{proposition}\label{prop:weak}
For every $t\in[0,1)$, the family $\mathcal W_t=\{[U]^>_t:U\in\Omega\}$ is a classical topology on $X$. The map $x\mapsto x_{1-t}$ is a homeomorphism of $(X,\mathcal W_t)$ onto the corresponding subspace of $(E,\tauQ)$.
\end{proposition}
\begin{proof}
The bottom and top cuts are $\varnothing$ and $X$. Strict inequality gives
\[
 [\bigvee_iU_i]^>_t=\bigcup_i[U_i]^>_t,\qquad
 [U\wedge V]^>_t=[U]^>_t\cap[V]^>_t.
\]
The homeomorphism follows because $x_{1-t}\in\Qq(U)$ if and only if $x\in[U]^>_t$.
\end{proof}

We now show that nonstrict cuts may fail to form a topology.

\begin{example} \label{ex:strongcuts}
Let $X=\mathbb N\cup\{\infty\}$. For $n\ge1$ define
\[
 U_n(k)=\begin{cases}1,&1\le k\le n,\\0,&k>n,\end{cases}
 \qquad U_n(\infty)=\tfrac12-\frac{1}{n+3}.
\]
Put $V=\bigvee_{n\ge1}U_n$, so $V(k)=1$ for $k\in\mathbb N$ and $V(\infty)=1/2$. The chain
\[
 \Omega=\{0,1,V\}\cup\{U_n:n\ge1\}
\]
is a Chang fuzzy topology. At level $1/2$, one has
\[
 [U_n]_{1/2}=\{1,\ldots,n\},\qquad [V]_{1/2}=X.
\]
Thus $\bigcup_n[U_n]_{1/2}=\mathbb N$ is not a nonstrict cut of an element of $\Omega$. The topology generated by nonstrict cuts therefore need not equal the family of nonstrict cuts itself. This explains why membership reconstruction uses $\tauM$ rather than assuming that $\Pp$ preserves arbitrary unions.
\end{example}

\subsection{Separation of fuzzy points}\label{sec:separation}
We close with two precise separation observations. They help distinguish the full point topologies from their fixed-level restrictions. A classical space is $T_0$ if any two distinct points are distinguished by an open set, and it is $T_1$ if each singleton is closed.

\begin{proposition} \label{prop:notT1}
If $X$ is nonempty, neither $(E,\tauM)$ nor $(E,\tauQ)$ is $T_1$.
\end{proposition}
\begin{proof}
Choose $x\in X$ and $0<\alpha<\beta\le1$. Every membership-open set containing $x_\beta$ contains $x_\alpha$, because all membership-open sets have downward closed vertical sections. Thus there is no membership-open set containing $x_\beta$ and excluding $x_\alpha$. Similarly, every quasi-open set containing $x_\alpha$ contains $x_\beta$, because all quasi-open vertical sections are upward closed. Both properties contradict the $T_1$ axiom.
\end{proof}

At a fixed level $\alpha$, let $\mathcal S_\alpha$ denote the classical topology on $X$ \emph{generated} by the nonstrict cuts $[U]_\alpha$, $U\in\Omega$, and let $\mathcal W_t$ be the strict-cut topology of Proposition~\ref{prop:weak}.

\begin{theorem} \label{thm:T0}
Suppose $\Omega$ is stratified. Then:
\begin{enumerate}[label=\textup{(\roman*)}]
\item $(E,\tauM)$ is $T_0$ if and only if $(X,\mathcal S_\alpha)$ is $T_0$ for every $\alpha\in(0,1]$;
\item $(E,\tauQ)$ is $T_0$ if and only if $(X,\mathcal W_{1-\alpha})$ is $T_0$ for every $\alpha\in(0,1]$.
\end{enumerate}
\end{theorem}
\begin{proof}
If either point space is $T_0$, every subspace $E_\alpha$ is $T_0$. The maps $x\mapsto x_\alpha$ identify the two restricted point topologies with $\mathcal S_\alpha$ and $\mathcal W_{1-\alpha}$, respectively, proving the necessary conditions.

For the converse, consider two distinct fuzzy points $x_\alpha$ and $y_\beta$. If $\alpha=\beta$ and $x\ne y$, the relevant fixed-level $T_0$ condition distinguishes them by an open set of the corresponding point space.

If $\alpha\ne\beta$, stratification supplies a constant fuzzy open set whose membership spectrum distinguishes the two levels: for example, if $\alpha<\beta$, the set $\Pp(\underline\alpha)$ contains $x_\alpha$ and excludes $y_\beta$. For quasi-coincidence, choose $c\in[0,1]$ with $1-\beta<c\le1-\alpha$; then $\Qq(\underline c)$ contains $y_\beta$ and excludes $x_\alpha$. Hence points at different levels are distinguished. This proves both sufficient conditions.
\end{proof}

\begin{remark}
Stratification is used only to distinguish unequal levels by constants; without it, the converse need not hold.
\end{remark}

\section*{Declarations}
\noindent\textbf{Funding.} No external funding was received for this work.\\
\textbf{Competing interests.} The author declares no competing interests.\\
\textbf{Data.} This is a theoretical study; no external dataset was used.

\noindent\textbf{Declaration of generative AI use.} OpenAI's ChatGPT was used to assist with language editing and LaTeX preparation The author takes responsibility for the final manuscript.

\end{document}